\documentclass[10pt,journal,twocolumn]{IEEEtran}
\usepackage{amsmath}
\usepackage[noend]{algorithmic}
\usepackage{algorithm}
\usepackage[final]{graphicx}
\usepackage{epstopdf}
\usepackage{epsfig, subfigure, bbm, amsfonts, latexsym, amssymb, url}
\usepackage{color}            % importation depuis Xfig

\usepackage[belowskip=-13pt,aboveskip=0pt]{caption}

\makeatletter
\newcommand*{\defeq}{\mathrel{\rlap{%
                     \raisebox{0.3ex}{$\m@th\cdot$}}%
                     \raisebox{-0.3ex}{$\m@th\cdot$}}%
                     =}
\makeatother

\newcommand{\argmax}{\operatornamewithlimits{argmax}}
\newtheorem{thm}{Theorem}
\newtheorem{lem}{Lemma}
\newtheorem{alg}{Algorithm}

\newcommand{\ignore}[1]{}

\newcommand{\crefrangeconjunction}{-}
\newcommand{\blue}[1]{{\color{blue} #1}}

\newif\iftp \tptrue

\title{\LARGE \bf
The Impact of Stealthy Attacks on Smart Grid Performance: Tradeoffs and Implications
}

\author{Yara Abdallah,~\IEEEmembership{Student Member,~IEEE,} Zizhan Zheng*,~\IEEEmembership{Member,~IEEE,} Ness B. Shroff,~\IEEEmembership{Fellow,~IEEE,} Hesham El Gamal,~\IEEEmembership{Fellow,~IEEE,} and Tarek M. El-Fouly,~\IEEEmembership{Member,~IEEE}% <-this % stops a space
\thanks{Y. Abdallah, Z. Zheng, and H. E. Gamal are with the Department of Electrical and Computer Engineering, The Ohio State University, 2015 Neil Ave., Columbus, OH 43210, USA. Email: yara.abdallah10@gmail.com; zheng.497@osu.edu; helgamal@ece.osu.edu.}
\thanks{N. B. Shroff is with the Department of Electrical and Computer Engineering, The Ohio State University, 2015 Neil Ave., Columbus, OH 43210, USA. Email: shroff.11@osu.edu.}
\thanks{T. M. El-Fouly is with the Department of Computer Science and Engineering, Qatar University, PO Box 2713, Doha, Qatar. Email: tarekfouly@qu.edu.qa.}
\thanks{* Corresponding author.}

\thanks{This work was supported by QNRF fund NPRP 5-559-2-227 and ARO-W911NF-15-1-0277. A preliminary version of this work appeared in the proceedings of the IEEE Conference on Decision and Control, 2012~\cite{Yara-cdc-12}.}
}

\begin{document}

\maketitle
%\thispagestyle{empty}
%\pagestyle{empty}

%%%%%%%%%%%%%%%%%%%%%%%%%%%%%%%%%%%%%%%%%%%%%%%%%%%%%%%%%%%%%%%%%%%%%%%%%%%%%%%%
\begin{abstract}
The smart grid is envisioned to significantly enhance the efficiency of energy consumption, by
utilizing two-way communication channels between consumers and operators. For example, operators
can opportunistically leverage the delay tolerance of energy demands in order to balance the energy load
over time, and hence, reduce the total operational cost. This opportunity, however, comes with security
threats, as the grid becomes more vulnerable to cyber-attacks. In this paper, we study the impact of
such malicious cyber-attacks on the energy efficiency of the grid in a simplified setup. More precisely,
we consider a simple model where the energy demands of the smart grid consumers are intercepted and
altered by an active attacker before they arrive at the operator, who is equipped with limited intrusion
detection capabilities. We formulate the resulting optimization problems faced by the operator and the
attacker and propose several scheduling and attack strategies for both parties. Interestingly,
our results show that, as opposed to facilitating cost reduction in the smart grid, increasing the delay tolerance of the
energy demands potentially allows the attacker to force increased costs on the system. This highlights
the need for carefully constructed and robust intrusion detection mechanisms at the operator.
\end{abstract}

\section{Introduction}\label{sec:introduction}

Over the past few years, the smart grid has received considerable momentum, exemplified in several regulatory and policy initiatives, and research efforts (see for example \cite{Moslehi2010,Lui2010} and the references therein). Such research efforts have addressed a wide range topics spanning energy generation, transportation and storage technologies, sensing, control and prediction, and cyber-security~\cite{McDaniel2009,false-data-CCS,data-attack-lang,malicious-price-lang,data-injection-Poor}.

Demand response/load balancing and energy storage are two promising directions for enhancing the energy efficiency and reliability in the smart grid. Non-emergency demand response has the potential of lowering real-time electricity prices and reducing the need for additional energy sources. The basic idea is that, by utilizing two-way communication channels, the \emph{emergency level} of each energy demand (at the end-users or central distribution stations) is sent to the grid operator that, in turn, \emph{schedules} these demands in a way that \emph{flattens} the load. Moreover, energy storage capabilities at the end-points offer more degrees of freedom to  the operator, allowing for a higher efficiency gain. This potential gain, however, comes at the expense of the security threat posed by the vulnerability of the communication channels to interception and impersonation.

In this work, we study the impact of the vulnerability of two-way communications on the energy efficiency of the smart grid. {More specifically, we propose a new type of data integrity attack towards Advanced Metering Infrastructures (AMI), that captures the above scenario in the presence of a single {\it stealthy} attacker. In an AMI system, a wide area network (WAN) connects utilities to a set of gateways, which are connected to electricity meters through neighborhood area networks (NANs). As observed in~\cite{AMI-threats-2012}, neighborhood area networks is an attractive target of attacks, where a large number of devices are physically accessible with little security monitoring available. Moreover, since these derives are connected to networks, an attacker can potentially get access to a large amount of data by hacking into a few nodes or links in AMI~\cite{AMI-threats-2010}. As observed in~\cite{Cleveland-AMI}, all the three major types of nodes in AMI, namely, smart meters, data concentrators, and the AMI headend, are subject to attacks, with different amount of data that can be utilized by the attacker.
}

In this work, we consider a simplified model of AMI, similar to \cite{Koutsopoulos2012}, that includes a grid operator and $n$ consumers that may be capable of energy storage, harnessing the potential cost savings in the smart grid. Our analysis covers two models of energy demands. In the first (total-energy model), each demand includes the total amount of energy to be served, the service start time, and the \emph{deadline} by which the requested energy should be delivered. In the second (constant-power model), each demand similarly has an arrival time and deadline, but the consumers ask for energy to be distributed across a specified number of time slots (a service time), with a power requirement in each slot. In both models, the consumers send their demands \ignore{simultaneously,}over separate communication channels to the operator. The grid operator attempts to schedule these demands so as to balance the load across a finite period of time, and hence \emph{minimize} the total cost paid to serve these demands.

\begin{figure}
 \centering
 \includegraphics[width=0.5\textwidth]{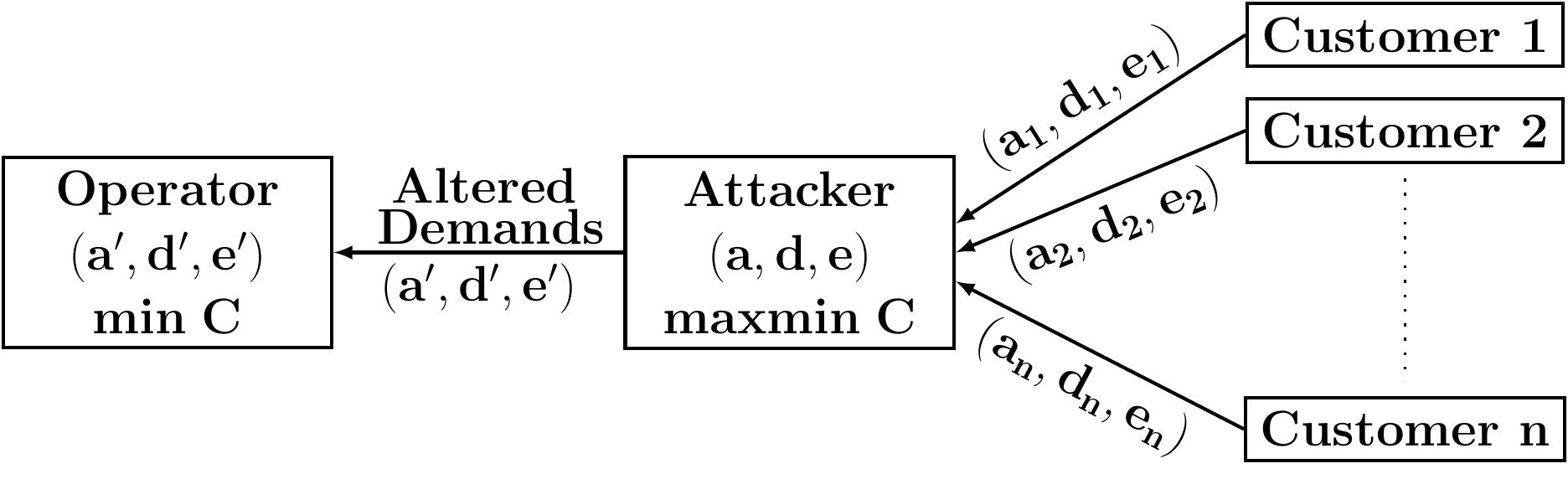}
 \caption{A system model for a smart grid in the presence of a single attacker. The forward channels between the consumers and the grid operator are fully compromised by the attacker. $(a,d,e)$ is the vector of the start times, deadlines and energy requirements of the consumer demands, respectively.}
 \label{fig:system_model}
\end{figure}

In our model, we also assume the presence of a single attacker who is fully capable of intercepting and altering the consumer demands before they arrive at the operator (see %, for example,
Figure~\ref{fig:system_model}). The end goal of the attacker, as opposed to the operator, is to \emph{maximize} the operational cost paid by the system for these demands, hence reducing the energy efficiency of the system. We differentiate between two scenarios. The first corresponds to a naive operator who fully trusts the incoming energy demands, whereas in the second, a simple intrusion detection mechanism (to be discussed later) is assumed to be deployed by the operator. Rather intuitively, the attacker's desire to remain undetected imposes more limitations on its capabilities, and hence, reduces the potential harm. This desire can be justified, for example, by considering the long-term performance of the grid, where the total impact of successive attacks is more damaging when the attacker remains undetected.

Based on the aforementioned assumptions, we first formulate the optimization problems faced by the operator and the attacker. %for both the total-energy and constant-power demand models.
For the operator, when being oblivious to any attacks, a minimization problem needs to be solved. On the other hand, the attacker is aware of the optimal strategy employed by the operator, and hence, a maximin optimization problem needs to be solved. In our formulation, we limit the attack's strength by the number of energy demands the attacker is capable of altering. %\emph{without} being detected.
For the case when the attacker is capable of altering \emph{all} of the energy demands (the attacks thus reach their full potential and force the system to operate at the maximum achievable total cost), we show that the maximin problem actually reduces to a maximization problem.

%Previous works on data integrality attacks in AMI mainly focus on energy theft and remote disconnection~\cite{AMI-anonmly,cybersecurity-smartgrid-survey}.
%Cyber attacks toward the smart grid have received increasing attention in recent years~\cite{McDaniel2009,Stuxnet,survey-1}. For instance, it has been observed that false data injection attacks can bias the power system state estimation~\cite{false-data-CCS,data-attack-lang}, leading to instable operation of the power grid as well as financial losses in electricity markets~\cite{malicious-price-lang,false-data-market}. Such attacks can be created if the attacker knows the configuration of the power system and can compromise a small fraction of smart meters. In the context of AMI, meter data and communication traffic, if not well protected, can disclose sensitive information and compromise the privacy of end users~\cite{McDaniel2009,smartmeter-privacy-Poor-1}. %Most of these works and other related works, however, mainly focus on the stable operation of the power grid.
To the best of our knowledge, however, the impact of stealthy attacks on the energy efficiency of the smart grid has not been studied before, and this paper is the first attempt to {\it explicitly} characterize such impact.

Our main contribution can be summarized as follows.

\begin{itemize}
    \item {We propose optimal offline strategies for both the operator and the unlimited attacker. The former gives the minimum energy cost when there is no attack, and the latter gives the maximum energy cost that can be enforced. % we propose optimal offline strategies. %where the optimization is carried over the entire set of the consumers' demands.
        The gap between the two indicates the maximum damage that can result from such attacks.} We also provide efficient online strategies for both of them. %where optimization is carried upon the arrival of each demand.
        These strategies are more practical in terms of operability and also indicate several bounds on the possible damage due to an unlimited attack. %An additional online attack is provided for the case when the attacker foresees future arrivals within a fixed lookahead window. This scenario captures the performance degradation when buffering capabilities are available to the attacker.
    \item For more limited attacks, we provide a simple greedy offline algorithm to arrive at a lower bound, and a dynamic programming-based algorithm that computes an upper bound on the total cost achieved by such attacks. Moreover, efficient online attacks are provided.  %with a characterization of their guaranteed performance.
    %\item We provide numerical results that support our theoretical claims under different scenarios. In these studies, we compare the average system performance in the presence/absence of attacks with the expected system performance when the delay tolerance of the jobs is not exploited by the operator (resembling the current electric gird where the communication infrastructure is absent). Moreover, we show the trade-off between the strength of the intrusion detection at the operator and the reduction in the system's efficiency due to stealthy attacks.
    \item From our analysis and numerical results, we conclude that {\bf in the absence of security threats} an increase in the delay tolerance of the energy demands increases the energy efficiency of the system, as expected, since the operator is offered more scheduling opportunities. On the other hand, a somewhat surprising observation is that, with a limited defense mechanism at the operator, this increase offers a similar opportunity to the attacker to force costs {\bf even higher than those incurred by the regular grid}, transposing the purpose of the communication capabilities provided to the consumers.
\end{itemize}

%[[.. maybe in Disucussion/Conclusion section:
The proposed framework enjoys several merits. Our analysis throughout the sequel does not assume any specific structure/distribution on the consumer demands and hence the derived results encompass a wide range of realistic scenarios. The attack bounds provided here are based on worst-case analyses and so provide strong guarantees on the impacts of different attacks. The main limitation of this work is the rather weak detection/defense mechanism at the operator. Our purpose here is to explore the attacker's side and arrive at performance bounds that motivate stronger defenses at the operator/consumers.

The remainder of this paper is organized as follows. After a brief overview of related work in Section~\ref{sec:related}, we present our system model and the optimization problems at the operator and attacker sides in Section \ref{sec:formulation}. In Sections \ref{sec:te} and~\ref{sec:st}, we provide offline and online attacks for the total-energy model and the constant-power model, respectively. %present both offline and online scheduling and attack strategies, where both full attack and limited attack strategies are derived. %we first describe optimal and suboptimal/online scheduling strategies for the operator in Section~\ref{sec:te_min}. We then study full attack and limited attack strategies for this model in Sections~\ref{sec:te_max_full} and~\ref{sec:te_max_limited}, respectively. In each case, we consider both offline and online settings. %In Section~\ref{sec:st} we shift our focus to the constant-power model: Section~\ref{sec:st_min} covers several scheduling strategies at the operator and discusses their performance. Section~\ref{sec:st_max_full} and~\ref{sec:st_max_limited} detail how to extend the previously studied attacks to the constant-power model.
%Section~\ref{sec:st} details how to extend these results to the constant-power model.
Numerical results are given in Section \ref{sec:numerical}. We provide some suggestions to the operator in Section~\ref{sec:rec}, whereas our conclusions are given in Section \ref{sec:conclusion}. \iftp Discussion of the model, extensions of our solutions to time-dependent cost functions, and details of some of the algorithms are provided in the Appendix. \else Extensions and missing proofs are provided in our online technical report~\cite{technical-report}.\fi

{
\section{Related Work}\label{sec:related}
Cybersecurity is of critical importance to the secure and reliable operation of the smart grid, which is challenging to achieve due to the large scale and the decentralized nature of the grid, the heterogeneous requirements of the components, and the coupling of the cyber and physical systems.
Various types of cyber attacks targeting the availability, integrity, and confidentiality of the smart grid have been studied, and both prevention and detection techniques have been proposed~\cite{cybersecurity-smartgrid-survey,AMI-threats-2012}.

Data integrality attack is considered as an important threat to the smart grid~\cite{cybersecurity-smartgrid-survey}. In particular, false data injection towards the SCADA systems has received a lot of attention recently~\cite{false-data-CCS,data-attack-lang}. By injecting malicious data into a small set of controlled meters, this attack can bias the state estimation of the system while bypassing the bad data detection in the current SCADA systems. Since the seminal work of~\cite{false-data-CCS}, much effort has been devoted to the problem of finding the minimum number of meters to be controlled to ensure undetectability~\cite{data-attack-lang,data-injection-Poor}. {Although this sparsest unobservable data attack problem is NP-hard in general, a polynomial time solution is given in~\cite{security-index-2014} for the case when the network is fully measured.} Moreover, strategic defense techniques have been developed~\cite{data-injection-Poor,hidden-attack} and the impact of data injection attack on real-time electricity market has been considered~\cite{malicious-price-lang,false-data-market}. When the attacker does not have enough number of controlled meters, a generalized likelihood ratio test is proposed to detect attacks~\cite{data-attack-lang}. {In addition to data attacks, the sparsest unobservable attack problem has been studied in closely related power injection attacks~\cite{Zhao-CDC-2013}.}

In the context of AMI, various potential threats have been identified~\cite{Cleveland-AMI,AMI-threats-2010,cybersecurity-smartgrid-survey}, including integrity attacks for the purposes such as energy theft and remote disconnection. Different intrusion detection systems (IDS) have been considered, including specification based~\cite{AMI-threats-2010,AMI-speficiation} and anomaly based approaches~\cite{AMI-anonmly}. In particular, a set of data stream mining algorithms are evaluated and their feasibility for the different components in AMI is discussed in~\cite{AMI-anonmly}. %where due to the lack of real AMI transaction data, the KDD Cup 1999 data set for intrusion detection in computer networks is used in evaluation.
The information requirements for detecting various types of attacks in AMI are discussed in~\cite{AMI-threats-2012}. Although data integrity attacks are considered as a potential threat in AMI~\cite{cybersecurity-smartgrid-survey}, its impact on the energy efficiency of the system has not been considered before, and proper intrusion detection schemes for the new type of attack that we consider remain open. %In particular, a specification as in~\cite{AMI-speficiation} for time-elastic demands, and
} 
\section{Problem Formulation}\label{sec:formulation}

\subsection{Demand Model}
In this paper, we adopt the control and optimization framework first proposed in~\cite{Koutsopoulos2012} for the demand side of the smart grid. This framework assumes a central operator and $n$ energy consumers that send their energy service demands to the operator using perfect channels. Our model builds on this framework, adding to it a single active attacker. The attacker is capable of intercepting and altering the demand requests in order to maximize the total energy cost paid by the smart grid. %To simplify the description, we assume $n = n_0$ in the rest of the paper.
We assume a time-slotted system and a finite time horizon $[0,T]$, and consider two types of demands:

\begin{enumerate}
  \item Total-Energy requirements: each consumer has a total \emph{energy} requirement that needs to be served before some deadline elapses. This, for instance, captures the scenario of having consumers with energy storage capabilities. Here, the energy demand of the $j^{th}$ consumer, $1 \leq j \leq n$, is composed of the tuple $(a_j,d_j,e_j)$, where $a_j,d_j \in \mathbb{N}^+, d_j \geq a_j, e_j \in \mathbb{R}^+$, indicating that %where $a_j,d_j \in \mathbb{N}^+$ denote the demand's arrival time and deadline, respectively, and $e_j \in \mathbb{R}^+$ denotes the requested total energy by the $j^{th}$ consumer.
      {demand $j$ arrives at the beginning of time-slot $a_j$, and has to be served for a total amount of $e_j$, by the end of time-slot $d_j$.}
  \item Constant-Power requirements: each consumer has an instantaneous power requirement and specifies a service time duration to finish a given job~\footnote{We use demand and job interchangeably in the paper.}, before a deadline elapses as well. The energy demand of the $j^{th}$ consumer, $1 \leq j \leq n$, is composed of the tuple $(a_j,d_j,s_j,p_j)$, where $a_j$ and $d_j$ are defined as in the above, $s_j \in \mathbb{N}^+$ is the job's duration time and $p_j \in \mathbb{R}^+$ is the instantaneous power requirement for this job. {We note that in contrast to the total-energy model, the instantaneous power requirement $p_j$ cannot be changed by either the operator or the attacker}.
\end{enumerate}

{In both cases, we assume that the set of jobs can be scheduled preemptively, i.e., a job can be interrupted and resumed, %after it starts,
so long as the deadline and energy/power requirements are met.} %We collect the set of job indices in
Let $J = \{1,\ldots,n\}$, and the associated demands in the total-energy (constant-power) model in $\mathcal{J} = \{(a_1,d_1,e_1),\ldots,(a_n,d_n,e_n)\}$ ($\mathcal{J} = \{(a_1,d_1,s_1,p_1),\ldots,(a_n,d_n,s_n,p_n)\}$). The set of demands are sorted by their arrival times non-decreasingly. \ignore{unless otherwise stated.} Each energy demand in $\mathcal{J}$ is sent to the operator over a perfect channel that is fully intercepted by the attacker. Hence the attacker could substitute the actual demand set $\mathcal{J}$ by a forged one, $\mathcal{J}'$, before it is received by the operator. An example for the total-energy case is shown in Figure~\ref{fig:system_model}. Similar to $\mathcal{J}$, the forged set is $\mathcal{J}' = \{(a'_1,d'_1,e'_1),\ldots,(a'_m,d'_m,e'_m)\}$ for the total-energy model and, for the constant-power model, is $\mathcal{J}' = \{(a'_1,d'_1,s'_1,p'_1),\ldots,(a'_m,d'_m,s'_m,p'_m)\}$. \ignore{We also let the corresponding forged job indices be} Let $J' = \{1,\ldots,m\}$ denote the indices of forged jobs. We note that $m \geq n$ in general as will be explained later in this section. {For ease of notation, we define the vector $a = [a_1,\ldots,a_n]$, and define $d,e,s,p$ similarly for the original vectors, and define $a',d',e',s',p'$ for the corresponding forged vectors. For any job $j$, we define its {\it job allowance} to be $l_j = d_j-a_j+1$.} Let $l_{max} = \max_{j \in J} l_j$, $l_{min} = \min_{j \in J} l_j$. We similarly define $e_{max}, e_{min}, s_{max}, s_{min}, p_{max}, p_{min}$.

{
\subsection{Simple Intrusion Detection}
We put the following constraints on the attacker. First, when the attacker chooses a job to modify, he is limited to changing its arrival time or its deadline time, or, breaking the job into multiple separate jobs (that would appear to the scheduler as independent jobs), so long as the final schedule is {\it admissible}. That is, all of the original jobs are served exactly their energy requirement (or service time and power requirement) upon or after their arrival and before or upon their real deadlines. Note that the attacker could easily be detected by the consumers if the final schedule is {\bf not} admissible.

Moreover, we assume that the operator adopts a simple statistical testing based intrusion detection scheme. For example, consider a statistical testing on the slackness of jobs. The slackness of a job $j$, denoted as $x_j$, is defined as the maximum time elasticity when serving the job. Formally, $x_j = l_i-1$ for the total-energy model, and $x_j = l_j-s_j$ for the constant-power model. Assume that the slackness of demands are $i.i.d.$ samples of a known distribution with mean $\mu$ and variance $\sigma^2$. For a set of $n$ demands received, the operator determines if it has been modified by using, for example, the one sample z-test with statistic $z = \frac{\overline{x}-\mu}{\sigma}\sqrt{n}$, with a significance level $\alpha$, the probability threshold below which the operator decides the data has been modified.

Assume that the attacker knows (1) the distribution of demands, and (2) the statistical testing and $\alpha$ adopted by the operator. %it can choose to modify a subset of demands to avoid detection.
If the attacker also knows $\sum_j x_j$ for the set of demands in $J$, %and $x_{max} = \max_j x_j$,
it can find the maximum amount of job slackness that can be reduced for demands in $J$, while still passing the z-test on slackness. When this knowledge is not available as in the more realistic online setting (to be precisely defined below), the attacker cannot ensure undetectablility. However, it can choose to modify a small number of jobs to ensure a small probability of detection, which is still useful to the attacker. %introduces a constraint on the total slackness of the modified jobs.
We can similarly consider a statistical testing on the arrive times or other parameters of demands. Instead of working on a constraint that depends on the concrete statistical testing used, we consider a simple constraint on the fractional of energy demands that the attacker is capable of altering, %without being detected (or being detected with a low probability),
which can be derived from the statistical testing used. In addition to simplifying the optimization problems for the attacker, such a bound can also be interpreted as a resource constraint to the attacker. We will consider other types of constraints in our future work. %This threshold is known a priori to all parties.
Let $B = \lfloor \beta n \rfloor, \beta \in [0,1]$ denote the number of jobs that the attacker can modify. %A small $\beta$ ensures a small probability of being detected, which is still useful to the attacker.

We note that an accurate statistical modeling of electric demands with time elasticity is by itself a challenging problem especially when the demands are correlated, which provides further opportunity to the attacker. Although the operator can also consider more advanced intrusion detection schemes such as data mining based anomaly detection, the high dimension of the data stream (large number of demands with overlapping durations) is a big challenge to be addressed. %can be applied, .  %techniques can be used for intrusion detection, their performance for the new type of attack we consider is clear as there is no real AMI transaction data that is public available. %On the other hand, the operator can choose to use a model-free intrusion detection approach, e.g., an active learning Although model-free approaches such as
}

\subsection{Optimization at the Operator and the Attacker}
Upon receiving the $m$ (altered) demands, $\mathcal{J}'$, an admissible schedule of these demands (jobs) is to be determined by the operator. %, where a schedule is admissible if each job is served its requested energy (or its requested service time slots with the exact power requirement in each time slot) upon or after its arrival time and before or upon its deadline. %\footnote{We put no limit on the total energy served in each time slot.}.
%Letting $T = \max_{j \in J'} d'_j$,
A schedule is given by $\mathcal{S} = [\mathcal{S}]_{jt} \in \mathbb{R^+}^m \times \mathbb{R^+}^T$, where $\mathcal{S}_{jt}$ denotes the amount of energy allocated to job $j$ in time-slot $t$. Let $E_{\mathcal{S}}(t)$ be the total energy consumed at time-slot $t \in [0,T]$ under schedule $\mathcal{S}$, i.e., $E_\mathcal{S}(t) = \sum_{j\in J'} {\mathcal{S}}_{jt}$. Let $C_t(E_\mathcal{S}(t))$ denote the cost incurred by the total power consumed at the time-slot $t$. We assume $C_t: \mathbb{R}^+ \rightarrow \mathbb{R}^+$ to be a general non-decreasing and convex function, as in~\cite{Koutsopoulos2012}. The convexity assumption resembles the fact that, as the demand increases, the differential cost at the operator increases, i.e., serving each additional unit of energy to increasing demand becomes more expensive~\cite{Koutsopoulos2012}. In our analysis and evaluation, we will consider the following commonly adopted power function as an example, where $C_t(E)=E^b, b\in \mathbb{R}, b\geq 1$, which allows for estimating the performance for a wide range of monotone increasing and convex functions. Moreover, for simplicity of exposition, we assume $C_t(\cdot)$ to be time invariant in the following and omit the subscript $t$. We show that most of our algorithms and analytic results can be extended to time-dependent cost functions in \iftp Appendix~\ref{sec:time-varing-cost} \else our technical report~\cite{technical-report}\fi. %Appendix~\ref{sec:time-varing-cost}.

The operator attempts to balance the load by finding an admissible schedule (given the altered demands by the attacker) that minimizes the total cost over the interval $[0,T]$. The optimization problem at the operator side, for the total-energy model, is then defined as follows:
%\begin{equation}
\begin{align*}
C_{min}(a',d',e') & =  \min_\mathcal{S}%\underset{\mathbf{x}}{\text{minimize}}
 \sum_{t=1}^{T} C(E_\mathcal{S}(t)) \tag{PminE}\label{pr:min} \\
 \text{s.t.} \hspace{1em}
&  \mathcal{S}_{jt} \geq 0, & \forall j \in J', \forall t \in [0,T],\\
&  \sum_{t=a'_j}^{d'_j} \mathcal{S}_{jt} = e'_j,  & \forall j \in J'.
\end{align*}
%\tag{PminE}\label{pr:min}
%\end{equation}
where we have dropped the constraint that no energy is served to a job $j$ outside $[a_j,d_j]$ since $C(\cdot)$ is monotone increasing. Similarly, the problem for the constant-power model is
%\begin{equation}
\begin{align*}
C_{min}(a',d',s',p') & =  \min_\mathcal{S}%\underset{\mathbf{x}}{\text{minimize}}
 \sum_{t=1}^{T} C(E_\mathcal{S}(t)) \tag{PminS}\label{pr:min2} \\
 \text{s.t.} \hspace{1em}
&  \mathcal{S}_{jt} \in \{0,p'_j\}, & \hspace{-3ex} \forall j \in J', \forall t \in [0,T],\\
&  \sum_{t=a'_j}^{d'_j} \mathbf{1}_{\mathcal{S}_{jt} = p'_j} = s'_j, & \hspace{-3ex} \forall j \in J'.
\end{align*}
%\tag{PminS}\label{pr:min2}
%\end{equation}
\noindent where $\mathbf{1}_{\mathcal{S}_{jt} = p'_j} = 1$ if $\mathcal{S}_{jt} = p'_j$, and is 0 otherwise. {The constraints in the both problems ensure the admissibility of the considered schedules.}

\ignore{
On the other hand, the attacker attempts to find appropriate values of $a',d',e'$ (or $a',d',s',p'$) in $\mathcal{J}'$ such that the cost achieved by the legitimate scheduler (at the operator side) is maximized, \emph{without} being detected. The intrusion detection capability at the operator is modeled as the number of energy demands the attacker is capable of altering without being detected (this is obtained, for instance, from statistical studies on previous demands). Note that this bound can also be interpreted as a resource constraint to the attacker. This
threshold is known a priori to all parties. Moreover, the attacker could easily
be detected at the consumers side if the final schedule %determined by the operator
is {\bf not} admissible.

These conditions form the constraints on the optimization problems at the attacker. When the
attacker chooses a job to modify, he is limited to changing its arrival time or its deadline time, or,
breaking the job into multiple separate jobs (that would appear to the scheduler as independent
jobs), so long as the final schedule is admissible. That is, all of the original jobs are served
exactly their energy requirement (or service time and power requirement) upon or after their
arrival and before or upon their deadline. In mathematical terms, we let $b_j$ be the collection of the (sub)jobs that the attacker generates out of job $j$, $1 \leq j \leq n$. Each (sub)job is, again, a tuple of the form $(a',d',e')$ or $(a',d',s',p')$.
}
On the other hand, the attacker attempts to find appropriate values of $a',d',e'$ (or $a',d',s',p'$) in $\mathcal{J}'$ such that the cost achieved by the %legitimate scheduler (at the operator side)
operator is maximized, subject to the number of demands that can be modified. %\emph{without} being detected. %In mathematical terms, we
Let $b_j$ be the collection of the (sub)jobs that the attacker generates out of job $j$, $1 \leq j \leq n$. Each (sub)job is, again, a tuple of the form $(a',d',e')$ or $(a',d',s',p')$. To guarantee an admissible final schedule, %as explained above,
each set $b_j$ should satisfy the following conditions:

In the total-energy model, for each job $j$:

\begin{subequations}\label{eq:te_const}
For $1 \leq k \leq |b_j|$
\begin{equation}\label{eq:te_const01}
a'_k, d'_k \in \mathbb{N}^+, e'_k \geq 0,
\end{equation}
\begin{equation}\label{eq:te_const02}
a_j \leq a'_k \leq d'_k \leq d_j,
\end{equation}
%and
\begin{equation}\label{eq:te_const03}
\sum_{1 \leq k \leq |b_j|} e'_k = e_j.
\end{equation}
\end{subequations}

In the constant-power model, for each job $j$:

\begin{subequations}\label{eq:st_const}
For $1 \leq k \leq |b_j|$
\begin{equation}\label{eq:st_const01}
a'_k, d'_k, s'_k \in \mathbb{N}^+,
\end{equation}
\begin{equation}\label{eq:st_const02}
a_j \leq a'_k \leq d'_k \leq d_j, p'_k = p_j,
\end{equation}
%and finally
\begin{equation}\label{eq:st_const04}
\sum_{1 \leq k \leq |b_j|} s'_k = s_j.
\end{equation}
\begin{equation}\label{eq:st_const05}
[a'_k,d'_k] \cap [a'_l,d'_l] = \emptyset, \quad \forall k,l, k \neq l, 1 \leq k,l \leq |b_j|.
\end{equation}
\end{subequations}

The sets $b_j$ are then collected in the forged demand vector, i.e., $\mathcal{J}' \colon = \bigcup_{1 \leq j \leq n} b_j$. Under this setting, the attacker solves the following optimization problems. For the total-energy model:
\begin{equation}
\begin{aligned}
C_{maxmin}(a,d,e,\beta) & = \max_{a',d',e',J^*}%\underset{\mathbf{x}}{\text{minimize}}
 C_{min}(a',d',e') \\
 \text{s.t.} \hspace{2em}
&    \text{Eqs \eqref{eq:te_const01} - \eqref{eq:te_const03}},\\
&    |J^*| \leq \beta n,
\end{aligned}
\tag{PmaxminE}\label{pr:maxmin}
\end{equation}
\noindent where $\beta \in \mathbb{R}, 0 \leq \beta \leq 1$, and \begin{equation}\label{eq:maximin_constraint}
   J^* = \left\{j\in J \colon b_j \neq \{(a_j,d_j,e_j)\} \right\}.
\end{equation}

\noindent Here $J^*$ denotes the set of consumer job indices that were \textbf{modified} by the attacker. %and $\beta$ denotes the {\it fraction} of jobs that can be altered without being detected, which reflects the intrusion detection capability of the system.
%The remainder of the constraints imply that, if the energy requirement of a job is not satisfied or a job is served outside its legitimate service duration, the attacker can be easily detected by the corresponding consumer. %\iftp{
In a similar fashion, we define the attacker's optimization problem for the constant-power model:
\begin{equation}
\begin{aligned}
C_{maxmin}(a,d,s,p,\beta) & = \max_{a',d',s',p',J^*}%\underset{\mathbf{x}}{\text{minimize}}
 C_{min}(a',d',s',p') \\
 \text{s.t.} \hspace{2em}
&    \text{Eqs (\ref{eq:st_const01}) - (\ref{eq:st_const05})},\\
&    |J^*| \leq \beta n,
\end{aligned}
\tag{PmaxminS}\label{pr:maxmin2}
\end{equation}
\noindent where $\beta \in \mathbb{R}, 0 \leq \beta \leq 1$, and
\begin{equation}\label{eq:maximin_constraint_service}
    J^* = \left\{j \in J \colon b_j \neq \{(a_j,d_j,s_j,p_j)\} \right\}.
\end{equation}
%}
%\else
%The attacker's optimization problem for the constant-power model can be defined similarly (omitted).
%\fi

We provide efficient {\it offline} and {\it online} solutions to the problems formulated above. Offline solutions not only give us performance bounds on the extreme case when there is no uncertainty on energy demands, but also provide useful insights for the design of online solutions. On the other hand, in the more realistic online setting, a demand is revealed only on its {\it actual} arrival. %Formally, we consider the following three cases.
{
\begin{itemize}
{\item Offline setting}: In the offline setting, we assume that the attacker knows all the true demands $\mathcal{J}$ at time 0, while the operator knows all the forged demands $\mathcal{J}'$ at time 0, and obtains no further information during $[0,T]$.
{\item Online setting}: In the online setting, at any time $t$, the attacker only knows the set of true demands with $a_j \leq t$, while the operator only knows the set of unmodified demands with $a_j \leq t$, and the set of forged demands with $a'_j \leq t$. In addition, the number of demands $n$ is the common knowledge.
%{\item Online setting with look-ahead}: To reduce the uncertainty about the future, the operator may choose to incentivize customers to share their future requirements through the two-way communication channels, which, however, can again be intercepted by the attacker. To model this scenario, we consider an online setting with a lookahead window $L \geq 0$, where at any time $t$, the attacker possesses knowledge about the demands with $a_i \leq t+L$, and the operator possesses knowledge about the forged demands with $a' \leq t+L$. Note that when $L = 0$, we get the online setting without look-ahead discussed above.
\end{itemize}
Note that in the online setting, if $a_j' = a_j$, demand $j$ should be forwarded to the operator without delay. On the other hand, if $a'_j > a_j$, the attacker should hold demand $j$ until $a'_j$ so that the operator does not get extra information.
}

For comparison purposes, we also consider the following inelastic scheduling policy for the operator as a baseline strategy. In the total-energy model, this strategy serves each job its energy demand, entirely and immediately upon its arrival. The associated baseline cost, $C_{base}(a,d,e)$, can be found as:
\begin{equation}\label{eq:c_base_TE}
C_{base}(a,d,e) = \sum_{t\in [0,T]} C\Bigg( \sum_{j\in J:a_j = t} e_j \Bigg).
\end{equation}

The counterpart quantity in the constant-power model is:
\begin{equation}\label{eq:c_base_ST}
C_{base}(a,d,s,p) = \sum_{t\in [0,T]} C\Bigg( \sum_{j\in J:t \in [a_j, a_j+s_j-1] } p_j \Bigg).
\end{equation}

This strategy %(under both models)
represents the case when the delay tolerance of the jobs is not exploited. Therefore, we treat this quantity as the cost paid in the \emph{current regular gird}, where no two-way communication channels are established, and accordingly, the system is not vulnerable to the cyber-attacks discussed in this paper.

%Finally, the following definitions are used throughout the remainder of this paper.
%\begin{itemize}
%  \item For any consumer job $j\in J$, define its \blue{job allowance to be $l_j = d_j-a_j+1$}. Let $l_{max} = \max_{j \in J} l_j$, $l_{min} = \min_{j \in J} l_j$, $e_{max} = \max_{j \in J} e_j$ and $e_{min} = \min_{j \in J} e_j$, \blue{$s_{max} = \max_{j \in J} s_j$ and $s_{min} = \min_{j \in J} s_j$, $p_{max} = \max_{j \in J} p_j$ and $p_{min} = \min_{j \in J} p_j$.}
%  \item %Given a job index set $J$ with size $n$, denote the set of the endpoints of the job intervals by $X_J \defeq \{a_1,\ldots,a_n\} \cup \{d_1,\ldots,d_n\}=\{1,\ldots,q\}$.
%      For every pair $(k,l), k \leq l$, let $\mathcal{I}_J(k,l)$ be the set of all job indices whose intervals are entirely contained in $[k,l]$, that is, $\mathcal{I}_J(k,l) = \{j \in J \colon a_j \geq k, d_j \leq l\}$.
%\end{itemize}

As a first attempt towards understanding the impact of stealthy attacks on smart-grid demand-response, we have made several simplifications in this work. In \iftp Appendix~\ref{modeling}\else our technical report~\cite{technical-report}\fi, we provide a discussion on the rationale behind our model and outline several extensions including how to conduct the impact analysis under congested power systems.

\section{total-energy Demands: Scheduling and Attack Strategies}\label{sec:te}
In this section, we focus on the total-energy demand model. We first find the optimal scheduling strategy for the operator %(the solution to Problem~\eqref{pr:min})
in Section \ref{sec:te_min}. We next propose full attack strategies in Section \ref{sec:te_max_full} including both offline and online attacks. Finally, in Section \ref{sec:te_max_limited}, we propose limited attacks and study the impact of such attacks. {We note that the offline attacks we discuss below have a time complexity of $O(n^3)$. On the other hand, all the online attacks have a time complexity of $O(n)$ and are therefore more scalable to large systems.}

\subsection{Scheduling at the Operator}\label{sec:te_min}
The optimization problem at the operator~\eqref{pr:min} can be directly mapped to the ``minimum-energy CPU scheduling problem" studied in \cite{Yao1995}. Our discussion below is an adapted discrete-time version of the classical YDS algorithm~\cite{Yao1995}.

For every pair $(k,l), k \leq l$, let $\mathcal{I}_J(k,l)$ be the set of all job indices whose intervals are entirely contained in $[k,l]$, that is, $\mathcal{I}_J(k,l) = \{j \in J \colon a_j \geq k, d_j \leq l\}$. For the received (forged) demands $\mathcal{J}'$, define the \emph{energy intensity} on $\mathcal{I}_{J'}(k,l)$ to be
\begin{equation}\label{eq:intensity}
    g(\mathcal{I}_{J'}(k,l)) = \frac{\sum_{j \in \mathcal{I}_{J'}(k,l)} e'_j}{l-k+1},
\end{equation}
\noindent Note that if we only consider the set of jobs in $\mathcal{I}_{J'}(k,l)$, a schedule that serves $g(\mathcal{I}_{J'}(k,l))$ amount of electricity in each time slot in the interval $[k,l]$ minimizes the energy cost (assuming it is admissible). We further define $(k^*,l^*) = \argmax_{(k,l): k \leq l} g(\mathcal{I}_{J'}(k,l))$, that is, $\mathcal{I}_{J'}(k^*,l^*)$ is the set of jobs with the maximum energy intensity among all $\mathcal{I}_{J'}(k,l)$ for any $k,l$ with $k\leq l$.

It is shown in \cite{Yao1995} that, for strictly convex $C(\cdot)$, the optimal strategy schedules a total energy of $g(\mathcal{I}(k^*,l^*))$ in each time slot in $[k^*,l^*]$. That is, the interval with the maximum energy intensity must maintain this intensity in the optimal schedule. This also implies that no jobs out of $\mathcal{I}(k^*,l^*)$ are scheduled with those in $\mathcal{I}(k^*,l^*)$. Hence a greedy algorithm that searches for $\mathcal{I}(k^*,l^*)$, schedules the jobs in $\mathcal{I}(k^*,l^*)$ and then removes those jobs (and the corresponding interval) from the problem instance, can be used to solve Problem~\eqref{pr:min}. The corresponding algorithm is outlined below (see~\cite{Yao1995} for the details).

\begin{algorithm}
\caption{Offline Scheduling at the Operator}\label{alg:min_offline}
{\small
\begin{algorithmic}[1]
\WHILE{$J' \neq \emptyset$}
\STATE $\mathcal{I}_{J'}(k^*,l^*) \gets$ an interval with the highest energy intensity;
\STATE Schedule the jobs in $\mathcal{I}_{J'}(k^*,l^*)$ according to the Earliest Deadline First (EDF) policy, such that $E_\mathcal{S}(t) = g(\mathcal{I}_{J'}(k^*,l^*))$, for all $t \in [k^*,l^*]$; %(which is always feasible).
\STATE Delete the jobs in $\mathcal{I}_{J'}(k^*,l^*)$ from $J'$ and modify the problem to reflect the deletion of jobs.
\ENDWHILE
\end{algorithmic}}
\end{algorithm}

%\begin{alg}\label{alg:min_offline}
%Repeat the steps below until $J'$ is empty.
%\begin{enumerate}
%  \item Identify $\mathcal{I}^*(k^*,l^*)$. Schedule the jobs in $\mathcal{I}^*(k^*,l^*)$, such that $E_\mathcal{S}(t) = E(\mathcal{I}^*(k^*,l^*))$, for all $t \in [k^*,l^*]$, according to the %Earliest Deadline First (EDF) policy (which is always feasible).
%  \item Modify the problem to reflect the deletion of the jobs in $\mathcal{I}^*$:
%
%  For all jobs $j\in J' \setminus \mathcal{I}^*$, if $a'_j \geq k^*$, set $a'_j \leftarrow \max(k^*-1, a'_j - (l^*-k^*) -1)$.
%
%  Modify $d'_j$ similarly. Set $J'\leftarrow J' \setminus \mathcal{I}^*$.
%\end{enumerate}
%\end{alg}

%\begin{alg}\label{alg:min_online}
%\begin{enumerate}
%  \item For each job $j \in J'$, compute $p_j = \frac{e'_j}{d'_j-a'_j+1}$. Set $p_j(t) = p_j$ for all $t \in [a'_j,d'_j], p_j(t) = 0$ otherwise.
%  \item For each time slot $t$, set $E_\mathcal{S}(t) = \sum_j p_j(t)$.
%\end{enumerate}
%\end{alg}

The above algorithm arrives at the optimal schedule with complexity $O(n^3)$ since it suffices to consider intervals whose two endpoints are either arrival times or deadlines of some jobs. Let $C_{min}$ denote the optimal minimum cost achieved (when there is no attack). A simple online algorithm for Problem~\eqref{pr:min} was also given in \cite{Yao1995} (the Average Rate Heuristic, AVR). This online scheme distributes the energy requirement of each job \emph{evenly} on its service interval, ignoring further information on how the jobs intersect. The performance of this simple heuristic is studied in \cite{Yao1995} when the cost mapping is a power function, and the following bounds are proven: For $C(E)=E^b, b\in \mathbb{R}, b\geq 2$, this online heuristic achieves a total cost $\overline{C}_{min} \leq r_b C_{min}$, where $b^b \leq r_b \leq 2^{b-1} b^b$. {Since each demand is processed once, this algorithm has an $O(n)$ complexity.}

%\begin{algorithm}
%\caption{Online Scheduling at the Operator}\label{alg:min_online}
%{\small
%\begin{algorithmic}[1]
%\STATE For each job $j \in J'$, compute $p_j = \frac{e'_j}{d'_j-a'_j+1}$. Set $p_j(t) = p_j$ for all $t \in [a'_j,d'_j], p_j(t) = 0$ otherwise.
%\STATE For each time slot $t$, set $E_\mathcal{S}(t) = \sum_j p_j(t)$.
%\end{algorithmic}}
%\end{algorithm}

%\begin{thm}\label{prop:min_online_tightness}
%For $C(E)=E^b, b\in \mathbb{R}, b\geq 2$, Algorithm \ref{alg:min_online} has an approximation factor of $r_b,$ where $b^b \leq r_b \leq 2^{b-1} b^b$.
%\end{thm}

\subsection{Full Attack Strategies and Performance Bounds}\label{sec:te_max_full}
We now turn our attention to the attacker's selection of $\mathcal{J}'$. We note that the special case ($\beta = 1$) is of special interest to us, as it resembles a \emph{full }attack, i.e., the attacker is capable of modifying \emph{all} of the consumer demands (e.g., when there is no intrusion detection at the operator). %without being detected.
We first address this case. The more general attacks for $\beta <1$ will be considered in Section~\ref{sec:te_max_limited}.

\subsubsection{An Optimal Offline Full Attack}\label{sec:te_max_offline_full}

We first show that, in the case $\beta =1$, the Problem~\eqref{pr:maxmin} can be transformed into a \emph{maximization} problem. To see this, consider any undetectable strategy followed by the attacker such that, for each demand $(a_j,d_j,e_j) \in \mathcal{J}$, there exists exactly one corresponding forged demand, $(a'_j,d'_j,e'_j) \in \mathcal{J'}$, with $a'_j = d'_j = t_j$ for some $t_j \in [a_j,d_j]$, and $e'_j = e_j$. All such strategies are always feasible to the attacker by our assumption of $\beta=1$ and, if employed by the attacker, leave no degrees of freedom to the operator. Moreover, due to the monotonicity and convexity of $C(\cdot)$, it suffices for the attacker to consider only this set of strategies as shown in the following lemma\iftp.\else~(see~\cite{technical-report} for the proof).\fi

\begin{lem}\label{lem:compress}
When $\beta =1$, there is an optimal attack where for any job $j$, $a'_j = d'_j = t_j$ for some $t_j \in [a_j,d_j]$.
\end{lem}
\iftp
\begin{IEEEproof}
Consider an optimal solution for the attacker. Suppose a job $j$ is served at both time $t_1$ and $t_2$. Let $E_1$ and $E_2$ denote the total energy consumption at $t_1$ and $t_2$, respectively. Without loss of generality, assume $C'_{t_1}(E_1) \geq C'_{t_2}(E_2)$. Then the total amount of $j$ served at $t_2$, denoted as $\delta$, can be moved from $t_2$ to $t_1$ such that $C_{t_1}(E_1+\delta) + C_{t_2}(E_2-\delta) \geq C_{t_1}(E_1) + C_{t_2}(E_2)$ by the convexity and monotonicity of $C_{t_1}$ and $C_{t_2}$. The lemma then follows by applying the above argument iteratively.
\end{IEEEproof}
\fi
%The above lemma holds even when the cost function $C(\cdot)$ varies over time as we prove in Appendix~\ref{proof_compress}.
Based on this observation, Problem~\eqref{pr:maxmin} under $\beta=1$ reduces to a maximization problem, which, for a given job instance, looks for an optimal strategy that serves each job in a \emph{single} feasible time slot. Formally, the attacker solves the following problem:
\vspace{-2ex}
\begin{equation}
\begin{aligned}
C_{max}(a,d,e) & =  \max_{\mathcal{S}}%\underset{\mathbf{x}}{\text{minimize}}
 \sum_{t=1}^{T} C(E_\mathcal{S}(t)) \\
 \text{s.t.} \hspace{1em}
&  \mathcal{S}_{jt} = 0,  & \hspace{-13ex} \forall j\in J, \forall t \in [0,T], t\neq t_j,\\
&  \mathcal{S}_{jt_j} = e_j, t_j \in [a_j,d_j], & \hspace{-13ex} \forall j\in J.\\
\end{aligned}
\tag{Pmax}\label{pr:max}
\end{equation}

%We will form the graph theoretic version of this problem, which is useful for describing the optimal full attack strategy, and for studying the impact of more limited attacks.

%Let $G = (V,E)$ be the interval graph induced by the job demands in $\mathcal{J}$. Each vertex $v_j \in V$ corresponds to a job interval that is given by $[a_j,d_j]$. An edge is thrown between any two vertices iff the two corresponding job intervals intersect at one or more time slots \cite{Gross2006}. In the induced interval graph, a \emph{clique} is a subset of vertices $S \subseteq V$, such that every two vertices in $S$ are connected by an edge. %A \emph{maximal clique} (inclusion-wise) is a clique that cannot be extended by including one more adjacent vertex.

Hence, in the above formulation, the attacker needs to decide only on $t_j$ for each $j\in J$. Given a set of jobs $J$, define a \emph{clique} of $J$ as a subset of jobs in $J$ whose job intervals intersect with each other, and a {\it clique partition} of $J$ as a partitioning of set $J$ into disjoint subsets where each subset forms a clique of $J$. We then have the following observation\iftp. \else~(see~\cite{technical-report} for the proof).\fi

\begin{lem}\label{lem:clique_partition}
Each clique partition of $J$ corresponds to a feasible solution to Problem~\eqref{pr:max} and vice versa.
\end{lem}
\iftp
\begin{IEEEproof}
Consider any clique partition of $J$. For each clique in the partitioning, the set of jobs in the clique overlap with each other, and can be compressed to the same time slot (any time slot where all these job intervals intersect). We then obtain a feasible solution to (\ref{pr:max}). On the other hand, consider a feasible solution to (\ref{pr:max}). We can assume that each job is served in a single time slot by Lemma~\ref{lem:compress}. For any time-slot $t$ with at least one job served, let $K_t$ denote the set of jobs that are served at $t$. Then $K_t$ is a clique for any $t$, and the set of these cliques form a clique partition of $J$.
\end{IEEEproof}
\fi

Moreover, we observe that to find the optimal attack, it is sufficient to consider {\it locally maximal cliques} defined as follows. For any time slot $t$, let $K^t$ denote the set of jobs whose job interval contains $t$. A clique is called {\it locally maximal} if it equals $K^t$ for some $t$. The following result is key to derive the optimal attack\iftp:\else~(see~\cite{technical-report} for the proof):\fi

\begin{lem}\label{lem:locally_maximal}
There is an optimal clique partition solving \eqref{pr:max} that contains a locally maximal clique~\footnote{A similar fact is proved in \cite{Gijswijt2007}, where the authors consider clique partitioning so as to \emph{minimize} a submodular cost function on the cliques, and shows the existence of a (globally) maximal clique in the optimal partition. We introduce the notion of locally maximal clique so that our results can be extended to time-dependent cost functions as we discuss in \iftp Appendix~\ref{sec:time-varing-cost}\else~\cite{technical-report}\fi.}.
\end{lem}
\iftp
\begin{IEEEproof}
Consider an optimal clique partition, $K_1, ..., K_m$, that solves \eqref{pr:max}. Assume $K_i$ has the maximum cost among these cliques. If $K_i$ is not locally maximal, then for any time-slot $t$ where jobs in $K_i$ intersect, there is a job $j$ included in another clique, say $K_{i'}$, whose interval contains $t$. By moving $j$ from $K_{i'}$ to $K_{i}$, we get a new partitioning whose total cost can only increase by the convexity and monotonicity of $C(\cdot)$. Hence, $K_i$ can be made locally maximal without loss of optimality.
\end{IEEEproof}
\fi

Let $\overline{C}(k,l)$ be the maximum feasible cost that could be achieved by solely scheduling the jobs in $\mathcal{I}_J(k,l)$. Given any time-slot $z$ contained in $[k,l]$, let $K^z_{k,l}$ be  the locally maximal clique at $z$ for jobs restricted to $\mathcal{I}_J(k,l)$. %i.e., $K_{k,l,z}$ is a maximal clique contained in $\mathcal{I}_J(k,l)$.
We then have the following recursion.

\begin{thm}\label{thm:recursion}
\begin{equation}
    \overline{C}(k,l) = \max_{z\in[k,l]} \Bigg[ C \Bigg( \sum_{j\in K^z_{k,l}} e_j \Bigg) + \overline{C}(k,z-1) + \overline{C}(z+1,l) \Bigg]. \label{eq:rec01}
\end{equation}
\end{thm}
\begin{IEEEproof}
Consider the set of jobs in $\mathcal{I}_J(k,l))$. Lemma~\ref{lem:locally_maximal} implies that $\overline{C}(k,l)$ is achieved by a partitioning that contains a locally maximal clique for jobs in $\mathcal{I}_J(k,l))$. Each such clique \emph{separates} the optimization problem into two subproblems for smaller intervals. By searching over all the locally maximal cliques over the interval $[k,l]$, $\overline{C}(k,l)$ can be achieved.
\end{IEEEproof}

Accordingly, we can apply the dynamic programming algorithm in~\cite{Gijswijt2007} to our problem as in Algorithm~\ref{alg:max_offline} (a formal description appears in \iftp Appendix~\ref{alg_max_offline}\else \cite{technical-report}\fi). The optimal cost is then $\overline{C}(1,T)$, which is computed in the final step together with the optimal clique partition. From the obtained clique partition, one can easily compute a set of time slots, $t_j, j\in J$ and set $a'_j = d'_j = t_j$, solving Problem~\eqref{pr:max}. The obtained schedule leaves no degrees of freedom to the operator as, after the attacker's modifications, all jobs become virtually urgent to operator and must be scheduled immediately upon their arrival. It is also clear that, as the job allowance of jobs increases, the attacker is capable of forming larger cliques and hence imposing higher costs on the operator. When we study online attacks, one of our goals is to formalize this observation.

\begin{algorithm}
\caption{Offline Full Attack}\label{alg:max_offline}
{\small
\begin{algorithmic}[1]
\STATE Iterate over all intervals $[k,l], k\leq l, k,l\in [0,T]$, with increasing interval length.
\STATE In each iteration, compute $\overline{C}(k,l)$ using Eq.~\eqref{eq:rec01}, where the last two terms are obtained from previous iterations.
\end{algorithmic}}
\end{algorithm}

{The algorithm has $O(n^2)$ iterations since it suffices to consider intervals whose two endpoints are either arrival times or deadlines of some jobs, where in each iteration, it takes $O(n)$ time to find $\overline{C}(k,l)$. Therefore, the algorithm has a total complexity of $O(n^3)$.}

%\begin{alg}[\cite{Gijswijt2007}]\label{alg:max_offline}
%\begin{enumerate}
%  \item Iterate over all intervals $[k,l], k\leq l, k,l\in X_J$, with increasing interval length.
%  \item In each iteration, compute $\overline{C}(k,l)$ using Eq.~\eqref{eq:rec01}, where the last two terms are obtained from previous iterations.
%\end{enumerate}
%\end{alg}

\ignore{
\begin{alg}[\cite{Gijswijt2007}]\label{alg:max_offline}
For all $k\in X_J$, set the initial condition
  \begin{equation}
     \overline{C}(k,k) = C \Bigg( \sum_{j\in \mathcal{I}_J(k,k)} e_j \Bigg).
  \end{equation}
With increasing subproblem width $(l-k)$, apply the following dynamic program:
\begin{enumerate}
  \item Compute
  \vspace{-1ex}
  \begin{equation}\label{eq:rec01}
    \overline{C}(k,l) = \max_{z\in[k,l]} \bigg[ C \left( \sum_{j\in K_{k,l}^z} e_j \right) + \overline{C}(k,z-1) + \overline{C}(z+1,l) \bigg],
  \end{equation}
   \noindent with $z^*$ achieving (\ref{eq:rec01}).

  \item Update the clique partition
  \begin{equation}\label{eq:rec02}
    \mathcal{Q}(k,l) = \begin{cases}
    \emptyset, \qquad \mbox{if } \mathcal{I}_J(k,l) = \emptyset,\\
    \mathcal{Q}(k,z^*-1) \cup K_{k,l}^{z^*} \cup \mathcal{Q}(z^*+1,l), \quad \mbox{otherwise.}
    \end{cases}
  \end{equation}
\end{enumerate}
\end{alg}
}

\subsubsection{Online Full Attacks} %($\beta=1$)}
In this section, we investigate the case where the attacker processes the arriving jobs in an online fashion, where at any time-slot $t$, the attacker possesses knowledge about the demands that have arrived by $t$. %The variable $L$ represents the potential batching/buffering capabilities at the attacker.
{We propose a simple online attack where the jobs in $J$ are partitioned into cliques according to an EDF policy. The attacker maintains a set of active jobs, that is, the set of demands that have arrived but not scheduled yet. Let $A$ denote the set of active jobs. In any time-slot $t$, if $t$ is the deadline for a demand $j \in A$, then %starting from the earliest deadline,
all the active demands in $A$ are grouped in a single clique, by setting their arrival times and deadlines to $t$. These demands are then forwarded to the operator, and $A$ is set to the empty set. Note that the algorithm ensures that the operator only learns a demand $j$ at $a'_j$. The intuition of using an EDF policy is to delay the decision as far as possible so that more demands can be compressed together to generate a large clique.}  % and then removed from the problem instance:

%\begin{alg}\label{alg:max_online}
%Repeat until $J$ is empty:
%\begin{enumerate}
%  \item Find the earliest deadline, $\tilde{d} = \min_{j \in J} d_j$.
%  \item Set $N = \{j\in J \colon a_j \leq \tilde{d} + L \}$.
%  \item Find the optimal clique partition of the jobs in $N$ using Algorithm \ref{alg:max_offline}, and set their arrival times and deadlines accordingly.
%  \item Update $J \leftarrow J\setminus N$.
%\end{enumerate}
%\end{alg}

\begin{algorithm}
\caption{Online Full Attack}\label{alg:max_online}
{\small
{
$A \gets \emptyset$. In any time-slot $t$,}
\begin{algorithmic}[1]
{
\STATE $A \gets A \cup \{j: a_j = t\}$;
\IF{$d_j = t$ for some job $j \in A$}
%\STATE $N \gets \{k\in J \colon a_k \leq d_j + L \}$;
%\STATE $N \gets \{k\in A \colon a_k \leq t\}$;
\STATE For each job $k$ in $A$, $a'_k \gets t, d'_k \gets t$;
\STATE Forward the set of (forged) jobs in $A$ to the operator;
%\STATE Find the optimal clique partition of the jobs in $N$ using Algorithm~\ref{alg:max_offline}, and set their arrival times and deadlines accordingly
\STATE $A \gets \emptyset$
\ENDIF
}
\end{algorithmic}}
\end{algorithm}

Since each job is processed once, the algorithm has a complexity of $O(n)$. We denote the resulting cliques by $K_1,\ldots,K_m$. The resulting cost is computed as
\begin{equation}\label{eq:max_bound}
    \underline{C}_{max} = \sum_{i=1}^m  C \bigg(\sum_{j \in K_i} e_j\bigg).
\end{equation}

Our next result shows that, despite its simplicity and online operation, Algorithm \ref{alg:max_online} could achieve a significant loss in the system's efficiency. We first make the following observation\iftp.\else~(see~\cite{technical-report} for the proof).\fi %This is shown in Theorem~\ref{prop:max_online_tightness} below (and our simulation results in Section~\ref{sec:numerical} as well).
\begin{lem}\label{lemma:max_online_tightness}
Consider any clique $X$ in an optimal (offline) solution that achieves $C_{max}$. Then $X = \bigcup_i (X \cap K_i)$, where $X \cap K_i$ and $X \cap K_j$ are disjoint for $i \neq j$, and $X \cap K_i$ is non-empty for at most $r_1$ different $K_i$, where $r_1 = \left\lceil{\frac{l_{\max}}{l_{\min}}}\right\rceil+1$.
\end{lem}
\iftp
\begin{IEEEproof}
Let $K_1,\ldots,K_m$ denote the sequence of cliques constructed by Algorithm~\ref{alg:max_online}. Since $K_i$ and $K_j$ contain disjoint set of jobs, and the union of all $K_i$ is the entire set of jobs, we have $X = \bigcup_i (X \cap K_i)$, where $X \cap K_i$ and $X \cap K_j$ are disjoint for $i \neq j$. Moreover, Algorithm~\ref{alg:max_online} ensures a property that for all $i'>i$, all the jobs in $K_{i'}$ have arrived strictly later than the earliest deadline of the jobs in $K_i$. Let $t_1$ and $t_2$ denote the earliest arrival and earliest dealine, respectively, among the set of jobs in $X$. Then since all the jobs in $X$ intersect at $t_2$, $t_2-t_1 \leq l_{max}$. The above property then ensures that $X$ could have a nonempty intersection with at most $r_1 \triangleq \lceil \frac{l_{max}}{l_{min}}\rceil + 1$ %\emph{consecutive}
sets in the partitioning $\{K_i\}, i\in\{1,\ldots,m\}$.
\end{IEEEproof}
\fi

This observation leads to the following bound for the online attack\iftp.\else~(see~\cite{technical-report} for the proof).\fi
\begin{thm}\label{prop:max_online_tightness}
For $C(E)=E^b, b\in \mathbb{R}, b\geq 1$, %Algorithm \ref{alg:max_online} has an approximation factor $r$, where
\begin{equation}\label{eq:max_online_tightness}
   \underline{C}_{max} \geq \frac{1}{{r_1}^{b-1}}C_{max}, \text{ where } r_1 = \left\lceil{\frac{l_{\max}}{l_{\min}}}\right\rceil+1
\end{equation}
\end{thm}
\iftp
\begin{IEEEproof}
For a given problem instance, $J$, $a,d,e$, let the optimal partition of the jobs in $J$ be $X_1, X_2, \ldots, X_{m^*}$, such that
\begin{equation}
    C_{max}(a,d,e) = \sum_{z=1}^{m^*} \Bigg( \sum_{j\in X_z} e_j \Bigg)^b.
\end{equation}

%We first consider the case $\delta = 0$.
Let $K_1,\ldots,K_m$ denote the sequence of cliques constructed by the algorithm. %We have that, for all $i'>i$, all the jobs in $K_{i'}$ have arrived strictly later than the earliest deadline of the jobs in $K_i$. Consequently, each $X_z, z\in \{1,\ldots,m^*\}$ could have a nonempty intersection with at most $r_1 \triangleq \lceil \frac{l_{max}}{l_{min}}\rceil + 1$ \emph{consecutive} sets in the partition $\{K_i\}, i\in\{1,\ldots,m\}$.
For any $z \in \{1,\ldots,m^*\}$, let $N(z,i) = X_z \cap K_i$. From Lemma~\ref{lemma:max_online_tightness}, we have
\begin{eqnarray}\label{eq:pmi}
    C_{max}(a,d,e)  &=& \sum_{z=1}^{m^*} \Bigg( \sum_{i=1}^m \Bigg( \sum_{j\in N(z,i)} e_j \Bigg) \Bigg)^b \nonumber \\
                    &\overset{(a)}{\leq}& \sum_{z=1}^{m^*} r_1^{b-1} \sum_{i=1}^m  \Bigg( \sum_{j\in N(z,i)} e_j \Bigg)^b \nonumber \\
                    &=& r_1^{b-1} \sum_{i=1}^m \sum_{z=1}^{m^*} \Bigg( \sum_{j\in N(z,i)} e_j \Bigg)^b \nonumber \\
                    &\leq& r_1^{b-1} \underline{C}_{max}(a,d,e), \label{eq:pmi}
\end{eqnarray}
\noindent where (a) is obtained by the power mean inequality.
\end{IEEEproof}
\fi
%The proofs are found in Appendix~\ref{proof_max}.
When $C(.)$ is a power function of the form $C(E)=E^b, b\in \mathbb{R}, b\geq 1$, the simple structure of the online solution further delivers an explicit lower bound for the maximum achievable cost by the attacker\iftp:\else~(see~\cite{technical-report} for the proof):\fi%, $C_{max}$:

\begin{thm}\label{prop:max_lower_bound}
For $C(E)=E^b, b\in \mathbb{R}, b\geq 1$,
\begin{equation}
    C_{max}(a,d,e) \geq \left( \frac{l_{min}\sum_{j \in J} e_j}{2 l_{min} + a_n - a_1} \right)^b.
\end{equation}
\end{thm}
\iftp
\begin{IEEEproof}
Suppose the attacker follows Algorithm~\ref{alg:max_online}. %with $L = 0$.
Let $K_1,\ldots,K_m$ denote the set of cliques constructed by the algorithm. From Eq.~\eqref{eq:max_bound} and the power mean inequality, we have
\begin{equation}
    C_{max} \geq \underline{C}_{max} \geq \left(\frac {\sum_{j\in J} e_j}{m}\right)^b. \label{eq:proof_max_lower_bound}
\end{equation}

Consider any two consecutive cliques $K_i$ and $K_{i+1}$. Let $j$ denote a job with the earliest deadline in $K_i$. Then from the construction of the algorithm, %and when $L = 0$,
we have $a_k-a_j \geq l_{min}$ for any job $k \in K_{i+1}$. Moreover, $a_1$ must appear in $K_1$ and $a_n$ must appear in $K_m$. It follows that $m \leq \frac{a_n-a_1}{l_{min}}+2$. This bound, together with~\eqref{eq:proof_max_lower_bound}, completes the proof.
\end{IEEEproof}
\fi

%\begin{figure}
% \centering
% \includegraphics[width=0.5\textwidth]{max-lower-bound.pdf}
% \caption{A lower bound on $C_{max}$ plotted for various values of $n$ and $l_{min}$ under a quadratic cost function (i.e., $b = 2$). The average energy demand is 10 while the average inter-arrival time is 5.}
% \label{fig:max_lower_bound}
%\end{figure}

%The proof is in Appendix~\ref{proof_max_lower}.
The above result formalizes our intuition that the harm done by a cyber attack grows with the scheduling leverage given to the grid's operator. When all other parameters are fixed, we use this result to specifically estimate the growth of $C_{max}$ with $l_{min}$. For instance, in Figure~\ref{fig:max_lower_bound}, %\iftp\else~(see the Appendix)\fi,
we plot our bound versus an increasing $l_{min}$, fixing the average energy demand and the average inter-arrival time. In this instance, $C_{max}$ grows at least linearly with $l_{min}$, and the rate of growth increases as the sample size $n$ increases. More numerical results are reported in Section \ref{sec:numerical}.

%\iftp
\begin{figure}[!t]
 \centering
 \includegraphics[width=0.4\textwidth]{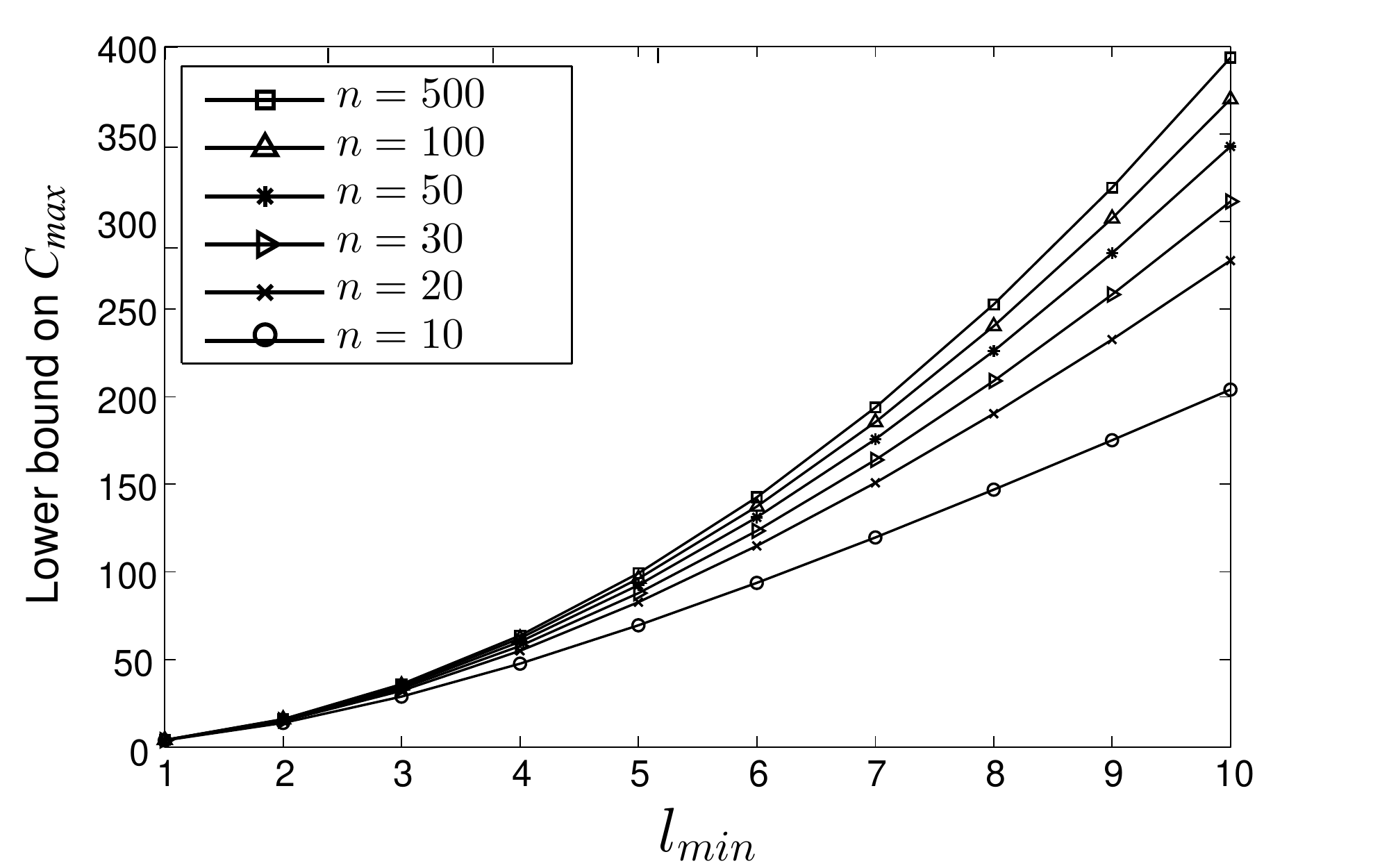}
 \caption{A lower bound on $C_{max}$ plotted for various values of $n$ and $l_{min}$ under a quadratic cost function (i.e., $b = 2$). The average energy demand is 10 while the average inter-arrival time is 5.}
 \label{fig:max_lower_bound}
\end{figure}
%\fi

\iftp
Finally, we use the heuristic AVR, presented previously in Section \ref{sec:te_min}, to arrive at an \emph{upper bound} on the gap between a fully-compromised operator (a operator subject to a \emph{full} attack) and a non-compromised one:

\begin{thm}\label{prop:max_upper_bound}
For $C(E)=E^b, b\in \mathbb{R}, b\geq 2$,
\begin{equation}
    C_{max}(a,d,e) \leq  2^{b-1} (l_{max}+1)^b   b^b C_{min}(a,d,e).
\end{equation}
\end{thm}
\begin{IEEEproof}
For a given problem instance, $J$, $a,d,e$, let the optimal partition of the jobs in $J$ be $X_1, X_2, \ldots, X_{m^*}$ such that
\begin{equation}
    C_{max}(a,d,e) = \sum_{z=1}^{m^*} \Bigg( \sum_{j\in X_z} e_j \Bigg)^b,
\end{equation}

\noindent and assume that those cliques are scheduled in time slots $t_1,\ldots,t_{m^*}$.

We now consider applying the AVR heuristic to the same problem instance and let the obtained cost by this algorithm be denoted by $C_{AVR}(a,d,e)$. We note that it achieves a fraction of $C_{max}(a,d,e)$ as given by the following:
\begin{eqnarray}
    C_{AVR}(a,d,e) &=& \sum_{t \in [0,T]} \Bigg( \sum_{j\in J} p_j(t) \Bigg)^b \\ \nonumber
    &\geq& \sum_{z=1}^{m^*} \Bigg( \sum_{j\in K_z} p_j(t_z) \Bigg)^b \\ \nonumber
    &=& \sum_{z=1}^{m^*} \Bigg( \sum_{j\in K_z} \frac{e_j}{l_j+1} \Bigg)^b \\ \nonumber
    &\geq&  \Bigg( \frac{1}{l_{max}+1} \Bigg)^b C_{max}(a,d,e).
\end{eqnarray}
In addition, we have $C_{AVR}(a,d,e) \leq 2^{b-1} b^b C_{min}(a,d,e)$~\cite{Yao1995} and that establishes our result.
\end{IEEEproof}
%The proof is provided in Appendix~\ref{proof_max_upper}.
This bound shows that, for a given power cost function,
the harm due to a cyber attack is also limited by the maximum scheduling flexibility given by
the served jobs.
\fi

\subsection{Limited Attacks and Performance Bounds}\label{sec:te_max_limited}
We now focus on the case when the attacker is limited by the number of jobs he is capable of modifying without being detected, i.e., the attacker can alter only $B = \lfloor \beta n \rfloor$ jobs, where $0 < \beta<1$. We again divide our study in two cases, the offline setting and the online setting. {In both cases, we derive bounds with respect to the corresponding full attacks. Therefore, these bounds are independent of the online scheduling algorithms used by the operator.}

\subsubsection{Offline Limited Attacks}\label{sec:te_offline_limited}
{For limited attacks, we are not able to find an optimal offline solution as we do for full attacks. To understand the impact of stealthy attacks in this more general setting, we propose two polynomial time offline algorithms that render a lower and an upper bound, respectively, on the performance of {\it optimal} limited attacks, and evaluate their performance in simulations. We show that even in the more challenging limited attack regime where the attacker may not be able to find the optimal attack, it is still possible to enforce significant amount of damage using a simple attack strategy.} %On the other hand, the upper bound the maximum damage from any attacker including those with unlimited computing power.}

Similar to our argument in the full attack case, the attacker could only consider the following simple strategy: Choose a set of job indices $J^* \subset J$ such that $|J^*|=\beta n$, and set $a'_j = d'_j = t^*_j$ for all the jobs $j\in J^*$. Leave all the remaining jobs $(J \backslash J^* )$ unaltered. We adopt this approach in our proposed offline attacks in this section. We let $C_{max} = C_{max}(a,d,e)$ and $C_{maxmin}(\beta) = C_{maxmin}(a,d,e,\beta)$, whenever clear from the context.

\iftp
\vspace{1ex}
\textbf{A Lower Bound.} \fi
{We first propose a simple variant that is tailored to our problem (see Algorithm~\ref{alg:maxmin_greedy}). For any $\beta$, the algorithm finds a feasible limited attack, the cost of which provides a lower bound on the cost resulting from the optimal limited attack. We further establish an explicit performance bound for this algorithm in Theorem~\ref{prop:maxmin_lower_bound}.}

Our algorithm is inspired by the standard greedy algorithm for the fractional knapsack problem~\cite{CLRS}. In the classical fractional knapsack problem, $m$ items are given, each with a weight $w_i$ and a value $v_i$. We need to find a set of items such that their total value is maximized subject to a budget on their total weight, say, $\beta_0 \sum_{i} w_i, 0\leq \beta_0 \leq 1$. A fraction of any item might be collected, and the corresponding value is scaled according to its chosen weight. The greedy algorithm below solves this problem.

%\begin{alg}\label{alg:greedy}
%Given the pairs $(v_1,w_1),\ldots,(v_m,w_m)$ and $\beta_0$
%\begin{enumerate}
%  \item Sort $(v_i,w_i)$ according to $v_i/w_i$ in a non-increasing order.
%  \item Choose the first $k$ pairs, $(v_1,w_1),\ldots,(v_k,w_k)$ such that
%  \begin{equation}
%      \sum_{i=1}^k w_i \leq \beta_0 \sum_{i=1}^m w_i \hspace{1em} \textrm{and} \hspace{1em}  \sum_{i=1}^{k+1} w_i > \beta_0 \sum_{i=1}^m w_i.
%  \end{equation}
%\end{enumerate}
%\end{alg}

%\noindent{Greedy knapsack algorithm:} Given the pairs $(v_1,w_1),\ldots,(v_m,w_m)$ and $\beta_0$
\begin{enumerate}
  \item Sort $(v_i,w_i)$ according to $v_i/w_i$ non-increasingly.
  \item Choose the first $k$ pairs, $(v_1,w_1),\ldots,(v_k,w_k)$ such that
  \begin{equation}
      \sum_{i=1}^k w_i \leq \beta_0 \sum_{i=1}^m w_i \hspace{1em} \textrm{and} \hspace{1em}  \sum_{i=1}^{k+1} w_i > \beta_0 \sum_{i=1}^m w_i.
  \end{equation}
\end{enumerate}

The optimal choice is given by the $k$ items collected in step (2), and a fraction of the $(k+1)^{th}$ item as the weight budget allows. Moreover, if we let the remaining weight budget after selecting the first $k$ pairs be parameterized by $\beta_1$ $= \left(\beta_0 \sum_{i=1}^m w_i - \sum_{i=1}^k w_i\right) / w_{k+1}$, by the greedy selection, we must have
\begin{equation}\label{eq:fractional-knapsack}
    \sum_{i=1}^k v_i + \beta_1 v_{k+1} \geq \beta_0 \sum_{i=1}^m v_i.
\end{equation}

The proposed attack strategy builds on the aforementioned algorithm:

\begin{algorithm}[!t]
\caption{Offline Limited Attack}\label{alg:maxmin_greedy}
{\small
\begin{algorithmic}[1]
{
\STATE Find the optimal clique partitioning using Algorithm \ref{alg:max_offline}, assuming a full budget;
\STATE Sort the set of cliques found by the ratio of clique cost over clique size non-increasingly;
\STATE Greedily choose a set of cliques with total size bounded by $\beta n$; Let $K$ denote the first unchosen clique on the list;
\STATE Greedily choose $\min(\beta n, |K|)$ jobs of highest energy requirements in $K$ to compress;
\STATE Among the above two choices, the one that results in a higher cost is adopted.
}
\end{algorithmic}}
\end{algorithm}

{Since finding the optimal clique partitioning is the most time-consuming step, the algorithm has a complexity of $O(n^3)$}. Let $C^1_{maxmin}(\beta)$ denote the total cost enforced by this attack. It is clear that $C^1_{maxmin}(\beta) \leq C_{maxmin}(\beta)$. To get insights on the performance of this attack, we consider two special cases. Suppose that, under no budget constraints, the optimal clique partition (obtained from Algorithm \ref{alg:max_offline}) is composed of cliques of size one, i.e., each job forms a separate clique. In this case, our greedy attack will choose to fully compress $B = \beta n$ jobs, and those will be of the highest energy demands according to step (3) above. By the greedy selection, this clearly guarantees that $C^1_{maxmin}(\beta) \geq \beta C_{max}$. Another extreme case is when the optimal clique partition is composed of one single clique containing all of the $n$ jobs. In this case, it is again clear by the greedy selection, in step (4), that $C^1_{maxmin}(\beta) \geq C\left(\beta \sum_{j\in J} e_j\right)$. When $C(.)$ is a power function of the form $C(E)=E^b, b\in \mathbb{R}, b\geq 1$, we get $C^1_{maxmin}(\beta) \geq \beta^b C_{max}$. For cases between those two extremes, we make use of the aforementioned insights to arrive at the following lower bound \iftp.\else~(see~\cite{technical-report} for the proof).\fi
%The algorithm provides us the following lower bound as proved in Appendix~\ref{proof_maxmin_lower}.

\begin{thm}\label{prop:maxmin_lower_bound}
For $\beta \in [0,1], C(E)=E^b, b\in \mathbb{R}, b\geq 1$,
\begin{equation}\label{eq:te_maxmin_bound_1}
    C^1_{maxmin}(\beta) \geq \frac{\beta^b}{2} C_{max}.
\end{equation}
\end{thm}
\iftp
\begin{IEEEproof}
Assume that the first $k$ cliques are fully compressed in Algorithm~\ref{alg:maxmin_greedy}. Let $\beta_1 = \frac{\beta n - (N_1+...+N_k)}{N_{k+1}}$ denote the fraction of budget available to clique $K_{k+1}$, where $N_i$ denotes the size of clique $K_i$. Let $C_0 = C_1 + \beta_1 E^b_{k+1}$. Then by the greedy selection of cliques and~\eqref{eq:fractional-knapsack}, we have $C_0 \geq \beta \sum_{i=1}^m E^b_i = \beta C_{max}$.

%\begin{equation}
%    C_0 \geq \beta \sum_{i=1}^m E^b_i = \beta C_{max}.
%\end{equation}

On the other hand, let $\beta_2 = \beta \frac{n}{N_{k+1}}$ denote the fraction of budget available to compressing only the jobs in clique $K_{k+1}$. By the greedy selection of jobs in the clique and~\eqref{eq:fractional-knapsack}, we have $\sum^{k'}_{j=1} e_j \geq \beta_2 E_{k+1}$. Therefore, $C_2 \geq \beta_2^b E^b_{k+1}$.
%\begin{equation}
%    C_2 \geq \beta_2^b E^b_{k+1}.
%\end{equation}

We then have
\begin{eqnarray*}
\frac{C^1_{maxmin}}{C_0} &=& \frac{\max(C_1,C_2)}{C_1 + \beta_1 E^b_{k+1}}
\geq \frac{C_2}{C_2 + \beta_1 E^b_{k+1}} \\
&\geq& \frac{\beta_2^b E^b_{k+1}}{\beta_2^b E^b_{k+1} + \beta_1 E^b_{k+1}}
= \frac{\beta^b_2}{\beta^b_2+\beta_1}  \\
&\overset{(a)}{\geq}& \frac{\beta^b_2}{\beta^b_2+\beta_2}
= \frac{\beta_2^{b-1}}{\beta_2^{b-1}+1}  \\
&\overset{(b)}{\geq}& \frac{\beta^{b-1}}{\beta^{b-1}+1}
\geq \frac{\beta^{b-1}}{2} ,
\end{eqnarray*}
\noindent where (a) follows from $\beta_1 \leq \beta_2$ and (b) follows from $\beta_2 \geq \beta$. Hence $C^1_{maxmin} \geq \frac{\beta^{b-1}}{2} C_0 \geq \frac{\beta^b}{2} C_{max}$.
\end{IEEEproof}
\fi

\iftp
\vspace{1ex}
\textbf{An Upper Bound.}
{In order to compute an upper bound on the maximum cost that can be obtained by {\it any} feasible offline limited attacks, %the system's performance,
we find the optimal attack strategy {\it under the assumption that the operator follows the baseline scheduling strategy}, i.e., the operator fully serves each job immediately upon its arrival. That is, we solve a problem similar to Problem \ref{pr:maxmin} by replacing $C_{min}$ with $C_{base}$. Let $C^2_{maxmin}$ denote the optimal total cost obtained by the attacker when the operator follows the baseline scheduling strategy. We first observe that $C^2_{maxmin}(\beta)$ is indeed an upper bound of $C_{maxmin}(\beta)$.}
{
\begin{lem}
$C^2_{maxmin}(\beta) \geq C_{maxmin}(\beta)$.
\end{lem}
\begin{IEEEproof}
Assume $S$ is the optimal attack strategy that achieves $C_{maxmin}(\beta)$, that is, $S$ is the optimal schedule for the attacker that solves Problem \ref{pr:maxmin}. Let $C'$ denote the total cost obtained when the attacker adopts $S$, while the operator adopts the baseline strategy. We then have $C_{maxmin}(\beta) \leq C' \leq C^2_{maxmin}(\beta)$.
\end{IEEEproof}
}

{
In the remainder of this section, we further assume that at most one job arrives at any given time-slot $t\in [0,T]$. Note that this is without loss of generality since we can consider an arbitrarily small time slot. We then show that under this assumption, $C^2_{maxmin}(\beta)$ can be found by a dynamic programming algorithm similar to Algorithm~\ref{alg:max_offline}.}
%The corresponding cost obtained (as an upper bound of $C_{maxmin}$) is evaluated in simulations in Section~\ref{sec:numerical}.

To derive the algorithm, we first observe that Lemma~\ref{lem:compress} and Lemma~\ref{lem:clique_partition} still hold for limited attacks, since the budget constraint is defined over the number of jobs that are altered, but not how they are altered. On the other hand, Lemma~\ref{lem:locally_maximal} does not hold any more. Instead, we will derive a variant of Lemma~\ref{lem:locally_maximal} as follows. Consider an optimal clique partition of $J$ to Problem~\eqref{pr:maxmin} when the operator follows the baseline scheduling strategy. %Let $O_1$ denote the set of all the cliques of size 1 in the partition and $O$ the rest of cliques.
Since at most one job can arrive at any time slot, without of loss of optimality, we can assume that each clique $K$ contains exactly one job, $j_K$, that has an unaltered arrival time. The remainder of the jobs would have arrival times altered to match that of $j_K$. For instance, we can choose $j_K$ as the job with latest arrival in clique $K$. Hence, the budget used to form clique $K$ would be exactly $|K|-1$. This observation leads to the following result. %(see~\cite{technical-report} for the proof). %Theorem \ref{prop:maxmin_DP} below. %\iftp(the proof is in Appendix~\ref{proof_maxmin_upper})\else(the proof is in our online technical report~\cite{technical-report})\fi.

\begin{thm}\label{prop:maxmin_DP}
When the operator follows the baseline strategy, there is an optimal clique partition solving Problem~\eqref{pr:maxmin} that contains a locally maximal clique, or a clique that can be made locally maximal by adding jobs from cliques of size 1 only.
\end{thm}
%\iftp
\begin{IEEEproof}
Let $K_{max}$ denote the clique containing the maximum total energy requirement in the clique partition. Assume that $K_{max}$ is not locally maximal. Then there exists a job $j$ contained in another clique $K$ in the partitioning such that $K_{max} \cup \{j\}$ is still a clique. Suppose $K$ contains at least 2 jobs. We distinguish the following two cases. First, if $a'_j \neq a_j$ in the optimal schedule, then we can schedule job $j$ at the time slot when all the jobs in $K_{max}$ are scheduled, while keeping the schedule of the rest of the jobs in $K$, without affecting the attacker's budget. Moreover, by the convexity of $C(.)$ and the fact that $K_{max}$ has the maximum total energy requirement among all the cliques in the partition, the resulting cost must increase by this change, which contradicts the fact that the clique partition is optimal. Second, if $a'_j = a_j$ in the optimal schedule, then we can again schedule job $j$ at the time slot when the jobs in $K_{max}$ are scheduled, and schedule the remaining jobs in $K$ at the latest arrival time of those jobs. By the assumption that at most one job arrives at any time slot and the fact that $|K| \geq 2$, this again leaves the budget unaffected and could only increase the total resulting cost. We again reach a contradiction. Hence, to achieve optimality in our upper-bound problem, what remains is to use jobs from cliques of size 1 to render $K_{max}$ maximal.
\end{IEEEproof}
%\fi

%The above theorem can be directly applied to any subgraph $G(\mathcal{I}(k,l))$, as defined in Section \ref{sec:te_max_offline_full}.
%Hence, similar to the case $\beta = 1$, for any such subgraph, each locally maximal clique contained in the subgraph separates the optimization problem into two subproblems.
Let $\overline{C}(k,l,m)$ denote the maximum achievable cost by solely scheduling the jobs contained in $[k,l]$ with a budget $m$. {Our objective is to find $\overline{C}(1,T,\lfloor \beta n \rfloor)$. Using Theorem~\ref{prop:maxmin_DP}, we can construct a recursion that computes $\overline{C}(k,l,m)$ by parsing for locally maximal cliques in each time-slot $z \in [k,l]$, as we did for $\beta = 1$, but with two modifications. First, we need to investigate all the possibilities of using only a fractional budget of $i$ out of $m$ for each found clique. Second, we would also need to exhaust the possibilities of distributing the remaining budget $m-i$ on the resulting two subproblems of any chosen clique.} Formally, for any clique $K$, let $K(i)$ denote the first $i$ jobs with the highest energy requirements in $K$. We then have:
\vspace{-1ex}
\begin{align}
    \overline{C}(k,l,m) &= \max_{z\in[k,l], i \in [0,m], j \in [0,m-i]} \Bigg[ C \Bigg( \sum_{j\in K^z_{k,l}(i+1)} e_j \Bigg) + \nonumber \\
    & \overline{C}(k,z-1,j) + \overline{C}(z+1,l,m-i-j) \Bigg]. \label{eq:rec02}
\end{align}

\noindent By Theorem~\ref{prop:maxmin_DP}, the constructed recursion indeed holds and a dynamic program similar to Algorithm \ref{alg:max_offline} is accordingly designed. {This algorithm has a complexity of $O(n^4)$, since it has $O(n^2)$ iterations and in each iteration it takes $O(n^2)$ time to find $\overline{C}(k,l,m)$.}
\else
{
We have further developed an algorithm that gives us an upper bound on the energy cost that can be achieved by any offline limited attacks. The details are given in our technical report~\cite{technical-report}.}
\fi

\subsubsection{Online Limited Attacks}
{To derive an efficient online limited attack, %we relax the budget constraint by allowing the attacker to adopt a random strategy so that the total number of modified demands is bounded by $\beta n$. Recall that we assume that the attacker knows $n$.
we consider the following simple strategy that mimics the behavior of Algorithm~\ref{alg:max_online} while taking the budget constraint into account. As in Algorithm~\ref{alg:max_online}, the attacker maintains the set of active jobs in $A$. It also maintains the total number of jobs that have been modified in $N$, and the number of future jobs in $R$ (recall that the attacker knows $n$). At any time $t$, the set of jobs that arrive at $t$ are added to $A$. The main idea of the algorithm is to modify each job with probability $\beta$, or forward it to the operator directly with probability $1-\beta$, independent of other jobs. Note that this decision has to made at the arrival time of a job. Let $A' \subseteq A$ denote the set of active jobs to be modified. If there is a job $j$ in $A$ with $d_j = t$, then all the jobs in $A'$ are compressed to the single time slot $t$. These jobs are then forwarded to the attacker, and both $A$ and $A'$ are set to the empty set. To make sure that all the budget is used and no more, the algorithm checks two boundary conditions. First, it stops sampling if all the budget has been used (lines 3-4). Second, when $R+N \leq B$, all the future jobs can be modified (line 6).} %Note that the algorithm ensures that the expected number of modified demands is bounded by $\beta n$.}

\begin{algorithm}[!t]
\caption{Online Limited Attack}\label{alg:maxmin_online}
{\small
{
$B \gets \lfloor \beta n \rfloor, A \gets \emptyset, A' \gets \emptyset, N \gets 0, R \gets n$. \\
In any time-slot $t$,}
\begin{algorithmic}[1]
{
\STATE $A \gets A \cup \{j: a_j = t\}$;
\FOR{each job $j$ with $a_j = t$}
\IF{$N=B$}
\STATE break;
\ENDIF
\STATE Sample $r$ from the uniform distribution in $[0,1]$;
\IF{$r \leq \beta$ or $R+N \leq B$}
\STATE $A' \gets A' \cup \{j\}$;
\ELSE
\STATE forward $j$ to the operator;
\ENDIF
\STATE $N \gets N + 1, R \gets R - 1$;
\ENDFOR
\IF{$d_j = t$ for some job $j \in A$}
%\STATE $N \gets \{k\in J \colon a_k \leq d_j + L \}$;
%\STATE $N \gets \{k\in A \colon a_k \leq t\}$;
\STATE For each job $k$ in $A'$, $a'_k \gets t, d'_k \gets t$;
\STATE Forward the set of forged jobs to the operator;
%\STATE Find the optimal clique partition of the jobs in $N$ using Algorithm~\ref{alg:max_offline}, and set their arrival times and deadlines accordingly
\STATE $A \gets \emptyset, A' \gets \emptyset$
\ENDIF
}
\end{algorithmic}}
\end{algorithm}

{
Since a separate decision is made for each demand on its arrival, and each demand to be modified is then processed once, this algorithm has a complexity of $O(n)$. Note that we have intentionally choose to generate the set of cliques at the earliest deadlines of jobs in $A$, not in $A'$, so that this algorithm closely simulates the behavior of Algorihtm~\ref{alg:max_online}. In particular, consider an input sequence, and any clique $K'$ generated by Algorithm~\ref{alg:maxmin_online}, and the corresponding clique $K$ generated by Algorithm~\ref{alg:max_online} at the same time slot. Then $K' \subseteq K$. Moreover, for a set of $i.i.d.$ demands, when $n$ becomes large, for most cliques $K$, the corresponding $K'$ has an expected size of $\beta |K|$. %  We have $K \subseteq K'$. Moreover, let $t_0$ denote the stoping time when the $\lfloor \beta n \rfloor$ jobs to be modified have been determined, or $t_0 = T$ if this does not happen in the given sample path. Then for every $K'$ generated before or upon $t$, we have E$(|K|) = \beta |K'|$.
Although there is no guarantee on the worst-case performance, we expect that the algorithm achieves an expected cost that is at least a constant fraction of $\left( \beta \frac{e_{min}}{e_{max}}\right)^b \underline{C}_{max}$ for $i.i.d.$ demands. %We note that this is not necessarily a tight bound as it is independent of the online cost minimization algorithm used by the operator.
}
\section{constant-power Demands: Scheduling and Attack Strategies}\label{sec:st}

Our previous scheduling and attack policies were derived solely for the total-energy demand model. In this section, we extend these results to demands that have service time and constant power requirements instead. We first provide an overview for the scheduling problem solutions at the operator in Section \ref{sec:st_min}. We then derive new full and limited attacks via simple modifications over the previously derived ones and analyze their performance in Sections \ref{sec:st_max_full} and \ref{sec:st_max_limited}, respectively.

\subsection{Scheduling at the Operator}\label{sec:st_min}
When all of the consumers require the same amount of power per time slot (i.e., $p_j = p$ for all $j \in J$), the Problem~\eqref{pr:min2} belongs to a class of ``load balancing" problems that are studied in detail in \cite{Hajek1990}. In this work, the author shows that the problem of finding the optimal schedule is equivalent to a network flow problem with convex cost. An optimal solution can be obtained by an iterative algorithm followed by a rounding step~\cite{Hajek1990}. For arbitrary power requirements, however, %even though the continuous version of this problem is solvable in polynomial time,
the integral nature of the problem renders it strongly NP-hard\iftp.\else~(see~\cite{technical-report} for a proof).\fi

\begin{thm}\label{prop:NP-hard}
For the constant-power model, Problem~\eqref{pr:min2} is strongly NP-hard.
\end{thm}
\iftp
\begin{IEEEproof}
We prove the result by a reduction from the 3-partition problem, which is known to be strongly NP-hard~\cite{Garey1979}. Consider an instance of the 3-partition problem: we are given a set $B$ of $3m$ elements $b_i \in Z^+, i = 1,...,3m$, and a bound $M \in Z^+$, such that $M/4 < b_i <M/2, \forall i$ and $\sum_i b_i = mM$. The problem is to decide if $B$ can be partitioned into $m$ disjoint sets $B_1,...,B_m$ such that $\sum_{b_i \in B_k} b_i = M$ for $1 \leq k \leq m$. Note that by the range of $b_i$'s, every such $B_k$ must contain exactly 3 elements. Given an instance of the 3-partition problem, we construct the following instance of our problem. There are $n=3m$ energy demands $J$, with $a_j = 1, d_j = m, s_j = 1$ and $p_j = b_j$ for all $j \in J$. The total power requirement of all consumers ($\sum_j p_j$) could be evenly distributed among the $m$ time slots if and only if the answer to the 3-partition problem is ``yes''. Clearly, such even distribution, if possible, corresponds to the optimal solution. Hence, solving Problem~\eqref{pr:min2} in this case answers the 3-partition problem,  making Problem \ref{pr:min2} strongly NP-hard.
\end{IEEEproof}
\fi

In our simulations, we report the relaxed continuous-version solution (as given in \cite{Hajek1990}) as a lower bound to the achieved cost by the optimal scheduler. {In this relaxed version, instead of a constant power $p_j$, job $j$ can be served by an amount $p_{jt} \in [0,p_j]$ for any time-slot $t$ such that $\sum_{t \in [a_j,d_j]} p_{jt} = s_jp_j$.} {We note that the continuous solution thus obtained can be furthered rounded to a feasible integral solution to the original problem. The main challenge, however, is to design the rounding process to achieve a low approximation factor, which remains open.}

As for online algorithms for the operator, solutions with performance guarantee are unknown for preemptive demands. Two scheduling policies were provided for non-preemptive demands in \cite{Koutsopoulos2012}. We choose the Controlled Release (CR) policy in our simulations, which is shown to be asymptotically optimal as average deadline duration approaches infinity \cite{Koutsopoulos2012}. In the CR policy, an active demand is served if the instantaneous power consumption in the current time slot is below a threshold or if it cannot be further delayed. Since each demand $j$ is processed at most $l_j$ times, independent of other demands, the algorithm has a complexity of $O(n)$. {Note that the online solution is always feasible and provides an upper bound to the offline optimal solution that is computationally hard to find.}

\subsection{Full Attack Strategies}\label{sec:st_max_full}

\subsubsection{Optimal Offline Full Attacks}\label{sec:st_offline_full}
In the case of full attacks ($\beta=1$), the total-energy demand model allowed the attacker to collapse the allowance of each job into a single time slot, while in this model, a job $j$ must be served in exactly $s_j$ time slots. However, we can still make use of the results developed earlier as follows. We break each job $j$ into $s_j$ separate sub-jobs, each having the same arrival time, deadline and the power requirement as those of $j$ and each should be served in exactly one time slot. With an entirely forced schedule on the operator, Problem~\eqref{pr:maxmin2} is thus turned into a maximization problem as before. In essence, to find the cost-maximizing schedule of those new (smaller) jobs, we are still attempting to form a clique partition of the resulting set of jobs %interval graph (with the same set cost function as before),
only with the additional constraint that no two subjobs resulting from a job $j$ can be scheduled in the same clique.

Let $\tilde{J} = \{(1,1),\ldots,(1,s_1),\ldots,(n,1),\ldots,(n,s_n) \}$ be the extended set of job indices, where $(j,k)$ denotes the $k$th subjob of the original job $j\in J$. Our clique partition is now over $\tilde{J}$. For any clique $K$, let $J_K = \{j \in J \colon (j,k) \in K, \text{ for some } k\}$, i.e., the set of jobs that originated the subjobs in $K$. For any time-slot $t$, we define a locally maximal clique, $K^t$, in this new setting as the set of subjobs that intersect at $t$, where at most one subjob from any job can be included. Following this definition, it is clear that the optimal solution indeed contains a locally maximal clique of subjobs, and, this also holds for any set of subjobs entirely contained within an interval.

{Let $\overline{C}(k,l,\{m_j\}_{j \in J})$ denote the maximum achievable cost by solely scheduling $m_j \leq s_j$ subjobs of job $j$ within interval $[k,l]$, which is defined to be 0 if for some $j$, $m_j > l-k+1$, or $m_j>0$ and $[k,l] \subsetneq [a_j,d_j]$. Our objective is to find $\overline{C}(1,T,\{s_j\}_{j \in J})$. Similar to Algorithm \ref{alg:max_offline}, we can construct a recursion that computes $\overline{C}(k,l,\{m_j\}_{j \in J})$ by parsing for locally maximal cliques in each time-slot $z \in [k,l]$. However, we observe that, unlike our previous model, a locally maximal clique in our extended set of jobs \emph{does not} divide a problem instance into a unique pair of smaller problems. Instead, all the potential subproblem-pairs resulting from a given locally maximal clique should be considered. We then have:
\vspace{-1ex}
\begin{align}
    \overline{C}(k,l,\{m_j\}_{j \in J}) &= \max_{z\in[k,l], m'_j \in [0,m_j-1] \forall j} \Bigg[ C \Bigg( \sum_{j\in K^z_{k,l}} p_j \Bigg) + \nonumber \\
    \overline{C}(k,z-1,&\{m'_j\}_{j \in J}) + \overline{C}(z+1,l,\{m_j-1-m'_j\}_{j \in J}) \Bigg]. \label{eq:rec03}
\end{align}
}

{We note that the complexity of this algorithm grows exponentially with the maximum clique size for a given problem instance, which indicates that the strategy can be computationally expensive for the attacker to use in practice. %However, this  does not rule out the possibility that a smart attacker can do better.
Due to the high complexity of the proposed attack, we have considered a relatively small scale setting in our simulations on offline attacks (see Figure~\ref{fig:bounds_var_powers}). % However, we find that, with large enough $n$, this parameter typically remains limited by the distribution parameters of arrival times and deadlines of the jobs.
An interesting open problem is to design a more efficient attack strategy that is close to optimal or rigorously prove that such an attack is hard to find.}

\subsubsection{Online Full Attacks}
{In the online case, we consider an attack similar to Algorithm \ref{alg:max_online}. The attacker again maintains a set of active jobs in $A$. In any time-slot $t$, the attacker checks if there is a job $j$ such that $d_j= t+s_j-1$. Note that to satisfy its service time requirement, such a job $j$ cannot be further delayed. If this is the case, all the jobs in $A$ are modified so that they will be scheduled for a consecutive number of time slots starting from $t$ until their service time requirements are satisfied. These jobs are then forwarded to the operator, and $A$ is the set to the empty set. It is important to notice that, similar to Algorithm \ref{alg:max_online}, if we only consider the set of jobs in $A$, then this strategy enforces the highest possible cost for those jobs.}%for served immediately is developed for the current demand model. The main difference is that to ensure the feasibility, the attacker needs to ensure that $d'$. Therefore, instead of the earliest deadline scheme must decide on the alteration of a given job $j$ upon $d_j-s_j+1$, since. The modified attack is then as follows:}

\begin{algorithm}
\caption{Online Full Attack (constant-power model)}\label{alg:max_online_st}
{\small
{
$A \gets \emptyset$. In any time-slot $t$,}
\begin{algorithmic}[1]
{
\STATE $A \gets A \cup \{j: a_j = t\}$;
\IF{$d_j= t+s_j-1$ for some job $j \in A$}
%\STATE $N \gets \{k\in J \colon a_k \leq d_j + L \}$;
%\STATE $N \gets \{k\in A \colon a_k \leq t\}$;
\STATE For each job $k$ in $A$, $a'_k \gets t, d'_k \gets t+s_k-1$;
\STATE Forward the set of (forged) jobs in $A$ to the operator;
%\STATE Find the optimal clique partition of the jobs in $N$ using Algorithm~\ref{alg:max_offline}, and set their arrival times and deadlines accordingly
\STATE $A \gets \emptyset$
\ENDIF
}
\end{algorithmic}}
\end{algorithm}

\ignore{
\begin{algorithm}
\caption{Online Full Attack (service time model)}\label{alg:max_online_st}
{\small
In any time-slot $t$,
\begin{algorithmic}[1]
%\STATE $J \gets J$ $\cup$ the set of jobs that arrive at $t$;
\STATE $N \gets \{k \in J \colon a_k \leq t \}$;
\IF{$d_j-s_j+1 = t$ for some job $j \in N$}
%\STATE $N \gets \{k \in J \colon a_k \leq d_j-s_j+1 + L \}$;
\STATE For each job $k$ in $N$, set $a_k \leftarrow \max (a_k,d_j-s_j+1)$;
\STATE Find the optimal clique partition of the jobs in $N$ for the constant-power model, and set their arrival times and deadlines accordingly
%\STATE $J \gets J \setminus N$
\ENDIF
\end{algorithmic}}
\end{algorithm}
}

%\begin{alg}\label{alg:max_online_st}
%Repeat until $J$ is empty:
%\begin{enumerate}
%  \item Set $\tilde{d} = \min_{j \in J} d_j-s_j+1$.
%  \item Set $N = \{j\in J \colon a_j \leq \tilde{d}+L \}$.
%  \item For each job $j$ in $N$, update $a_j \leftarrow \max (a_j,\tilde{d})$.
%  \item Find the optimal clique partition of the jobs in $N$ using the optimal offline full attack algorithm, and set their arrival times and deadlines accordingly.
%  \item Update $J \leftarrow J\setminus N$.
%\end{enumerate}
%\end{alg}

%With a lookahead window of $L = \delta l_{max}$, $0 \leq \delta \leq 1$, and
{
Since each demand is processed once, this algorithms has a complexity of $O(n)$. Similar to Lemma~\ref{lemma:max_online_tightness}, we have the following observation for the constant-power model.
\begin{lem}\label{lemma:max_online_partial_tightness}
Any clique $X$ in an optimal (offline) solution that achieves $C_{max}$ is a disjoint union of $X \cap K_i$, where $X \cap K_i$ is non-empty for at most $r_2$ different $K_i$, where $r_2 = (s_{max}-s_{min}+1)\left(\left\lceil{\frac{\max_j(l_j-s_j+1)}{\min_j(l_j-s_j+1)}}\right\rceil+1\right)$.% if $s_j = s$ for all $j$, and $r_2 \leq 2l_{max}$ in general.
\end{lem}

It is then straightforward to extend the proof of Theorem~\ref{prop:max_online_tightness} to show that the above algorithm achieves at least a fraction $\frac{1}{r_2^{b-1}}$ of the optimal offline cost in the constant-power model.
}

\subsection{Limited Attacks} \label{sec:st_max_limited}

\subsubsection{Offline Limited Attacks}
{To derive an offline limited attack in the constant-power model, we consider an algorithm similar to Algorithm \ref{alg:maxmin_greedy}. The optimal offline algorithm discussed in the previous section is first applied to find the optimal clique partitioning of sub-jobs when there is no budget constraint. Greedy algorithms are then applied twice; once to choose a set of cliques to fully compress, and to choose a set of sub-jobs within the first unchosen clique on the list, and the choice that results in a higher cost is adopted. In both cases, we require the total number of sub-jobs chosen to be bounded by $\beta n$. This ensures that the total number of modified jobs is also bounded by $\beta n$. Since finding the optimal clique partitioning may take exponential time in the worst case, this algorithm also has an exponential time complexity. Let $s_{avg} = (\sum_j s_j)/n$ denote the average service time requirement. Assume $C(E)=E^b, b\in \mathbb{R}, b\geq 1$. Since there are $\sum_j s_j$ sub-jobs in total and $\beta n$ of them are compressed, following Theorem \ref{prop:maxmin_lower_bound}, the guaranteed performance of this attack readily becomes $C^1_{maxmin}(\beta) \geq \frac{1}{2} \left(\frac{\beta}{s_{avg}}\right)^b C_{max}$.}

%To establish a guaranteed lower bound on limited attacks ($\beta<1$), an attack similar to Algorithm \ref{alg:maxmin_greedy} could be utilized for this demand model. We first note that $B = \beta n = \frac{\beta n}{\sum_j s_j} \times \sum_j s_j = \left(\frac{\beta}{s_{avg}}\right) \times \sum_j s_j,$ where $s_{avg}$ is the average service time requirement. In our attack, the optimal clique partition (over subjobs) is first computed. Accordingly cliques can be greedily chosen to be forced on the operator using Algorithm \ref{alg:maxmin_greedy}, with a budget fraction $\left(\frac{\beta}{s_{avg}}\right)$. For $C(E)=E^b, b\in \mathbb{R}, b\geq 1$, following Theorem \ref{prop:maxmin_lower_bound}, the guaranteed performance of this attack readily becomes $\underline{C}_{maxmin}(\beta) \geq \frac{1}{2} \left(\frac{\beta}{s_{avg}}\right)^b C_{max}$.

%For $C(E)=E^b, b\in \mathbb{R}, b\geq 1$,
%{\small
%\begin{equation}\label{eq:st_maxmin_bound}
%    \underline{C}_{maxmin}(\beta) \geq \frac{1}{2} \left(\frac{\beta}{s_{avg}}\right)^b C_{max}.
%\end{equation}
%}
%\end{thm}

\subsubsection{Online Limited Attacks}
{
We then modify Algorithm~\ref{alg:max_online_st} to obtain an online limited attack as we did for the total-energy model. The attacker maintains the set of active jobs in $A$, and samples a fraction $\beta$ of them to be modified, saved in $A'$. At any time $t$, if $d_j = t+s_j-1$ for some job $j$ in $A$, all the job in $A'$ are modified as in Algorithm~\ref{alg:max_online_st}. The algorithm also checks the two boundary conditions as we explained before to ensure that all the budget is used and no more. Assume $C(E)=E^b, b\in \mathbb{R}, b\geq 1$. Similar to Algorithm~\ref{alg:maxmin_online}, this algorithm also has a complexity of $O(n)$. As in the total-energy model, although there is no worst-case guarantee, we expect that this simple attack obtains an expected cost that is at least a constant fraction of $\underline{C}_{maxmin}(\beta) \geq \left( \beta \frac{p_{min}}{p_{max}}\right)^b \underline{C}_{max}$ for $i.i.d.$ demands and when $s_j$ is a constant for all $j$.} %\frac{s_{min}}{s_{max}}\left( \beta \frac{p_{min}}{p_{max}}\right)^b \underline{C}_{max}$.}

%Moreover, with a budget fraction $\beta$ and a lookahead window of $L = \delta l_{max}$, $0 \leq \delta \leq 1$, and for $C(E)=E^b, b\in \mathbb{R}, b\geq 1$, we get an approximation factor of $\left( \frac{\beta}{2} \frac{e_{min}}{e_{max}}\right)^b \left( \max \left(\frac{1}{2 l_{max}},\frac{\delta}{2}\right)\right)^{b-1}$, by extending Theorem~\ref{prop:maxmin_lower_bound_lookahead} to the constant-power model.

%\input{practical}
%%%%%%%%%%%%%%%%%%%%%%%%%%%%%%%%%%%%%%%%%%%%%%%%%%%%%%%%%%%%%%%%%%%%%
\section{Numerical Results}\label{sec:numerical}
In this section, we provide numerical results that illustrate the impact of stealthy attacks under various settings. \iftp\else More results can be found in our online technical report~\cite{technical-report}.\fi In this section, unless stated otherwise, the job arrivals are simulated as a Poisson arrival process with mean 3. We use a quadratic cost function $C(E) = E^2$ in all of our simulations.

\begin{figure}[!t]
\centering
\begin{tabular}{cc}
\begin{minipage}{0.47\linewidth}
  \centerline{\includegraphics[width=1\linewidth]{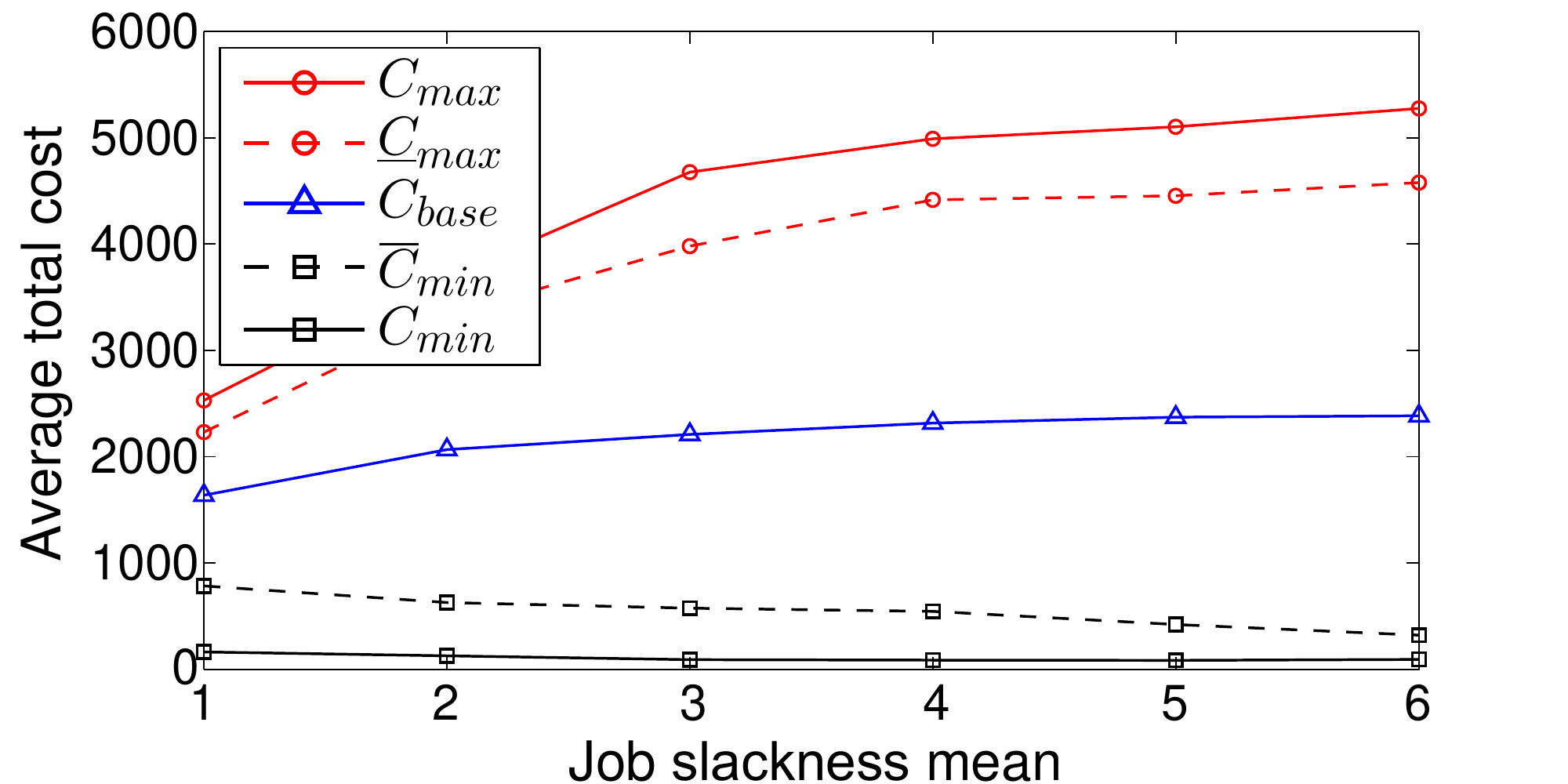}}
  \centerline{\small (a) total-energy model}
\end{minipage}
%\hfill
\begin{minipage}{0.47\linewidth}
  \centerline{\includegraphics[width=1\linewidth]{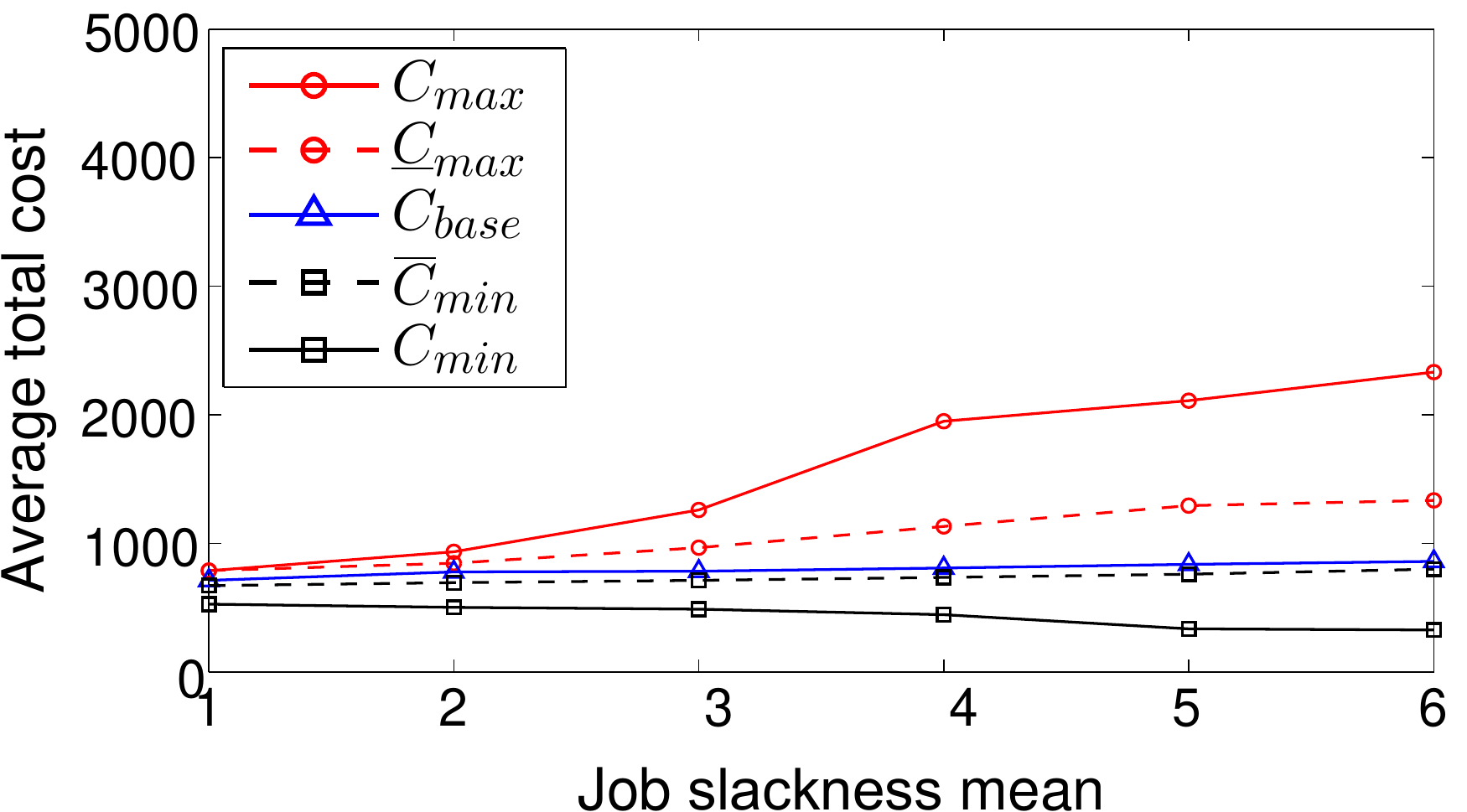}}
  \centerline{\small (b) constant-power model}
\end{minipage}
\end{tabular}
\caption{Comparison between the performance of a fully-compromised smart grid (offline and online attacks), the current grid, and an un-compromised smart grid (offline and online scheduling), under varying job allowance means.}
\label{fig:bounds_var_powers}
\end{figure}

\vspace{1ex}
\noindent{\bf Full Attacks:} In Figure~\ref{fig:bounds_var_powers}, we compare the performance of a non-compromised smart grid, a fully-compromised smart grid and the ``dumb'' grid (where all jobs are immediately scheduled upon their arrival), for both the total-energy model and the constant-power model, for a total of 20 jobs. All the job slackness are $i.i.d.$ exponential random variables, as well as the service time intervals. In the constant-power model, the job slackness mean is varied between 1 and 6, and the service time mean is fixed to 2. The power requirement per time slot, for each job, is uniformly distributed in the interval $[1,5]$. For comparison purpose, for each job generated in the constant-power model, a job with the same arrival, slackness, and total power requirement is generated for the total-energy model. The plots report the average performance of both systems over 10 trials. For the total-energy model, $C_{min}, \overline{C}_{min}, C_{base}, \underline{C}_{max}$, and $C_{max}$ correspond to the cost achieved by Algorithm~\ref{alg:min_offline}, the AVR algorithm, the baseline cost~\eqref{eq:c_base_TE}, Algorithm~\ref{alg:max_online}, and Algorithm~\ref{alg:max_offline}, respectively. For the constant-power model, they correspond to the lower bound obtained from the continuous relaxation of the minimization problem for the operator, the cost obtained by the Controlled Release (CR) policy~\cite{Koutsopoulos2012}, the baseline cost~\eqref{eq:c_base_ST}, the cost obtained by Algorithm~\ref{alg:max_online_st}, and that by the optimal offline full attacks discussed in Section~\ref{sec:st_offline_full}, respectively.

We observe that, as the job slackness mean increases, for both models, further scheduling opportunities are offered to the
legitimate operator, and hence further savings in the total cost are attained if the smart grid is not compromised. In the presence of an attacker, however, a similar flexibility is available to the attacker, and accordingly the severity of the attack increases as the job slackness mean increases. We also observe that the uncompromised total-energy system outperforms the constant-power model, in terms of total cost, due to the increased job scheduling flexibility in the former.
For the same reason, attacks are more harmful for this model as well. In the total-energy model, when compared to the costs paid by the regular grid, an offline (online) attack causes an increase in cost by 154\% (136\%) with a job slackness mean of 1 and up to 220\% (191\%), while the expected cost to be paid for an uncompromised system should, in fact, decrease by values ranging in $200\% - 2500\%$. A similar comparison could be drawn in the constant-power model. Therefore, overall, the unprotected smart grid simulated here, not only does it fail to meet the cost savings prospected in a smart grid, it performs far worse than the current electric grid.

\begin{figure}[!t]
\centering
\begin{tabular}{cc}
\begin{minipage}{0.47\linewidth}
  \centerline{\includegraphics[width=1\linewidth]{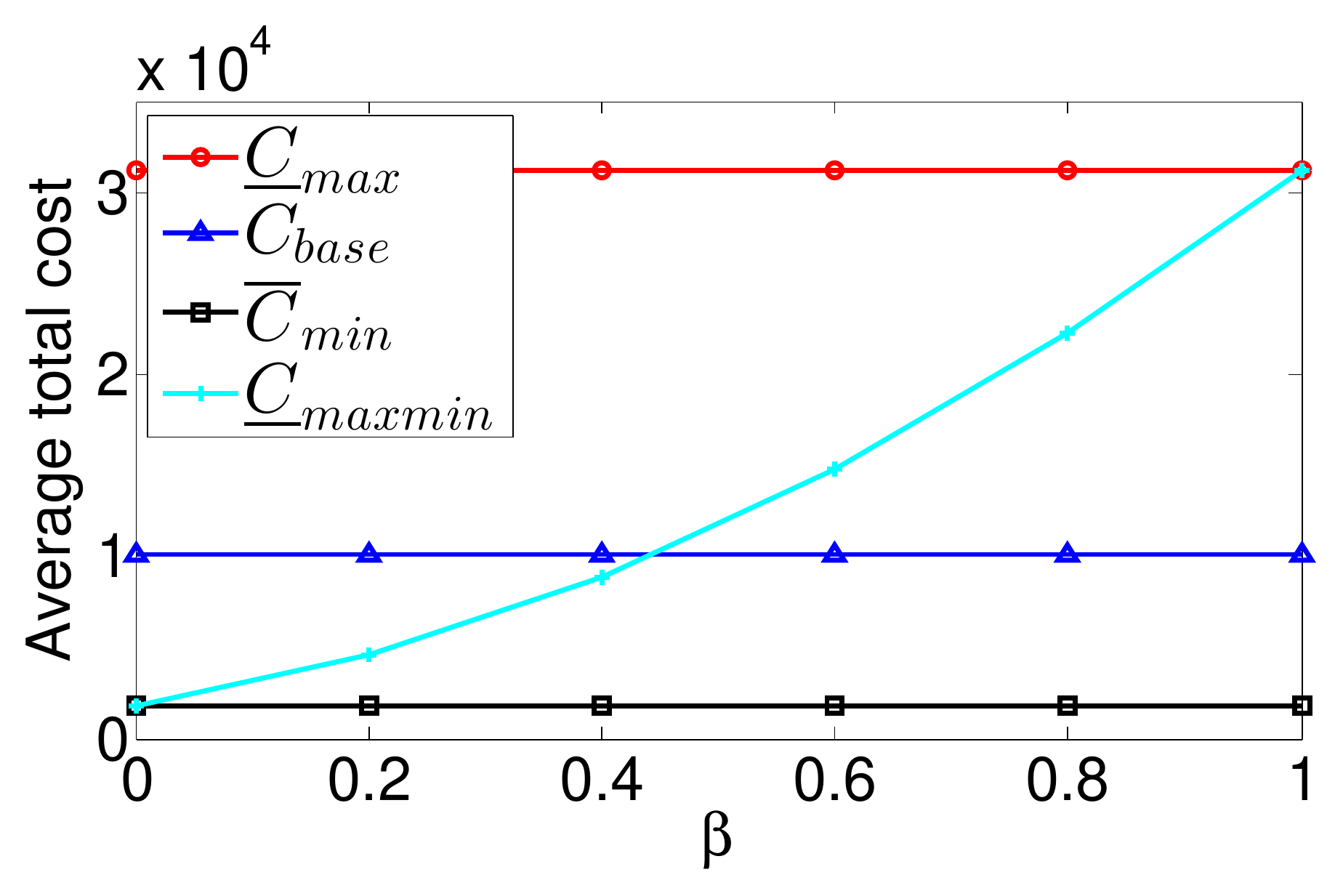}}
  \centerline{\footnotesize $x_i \in [0,40]$}
\end{minipage}
%\hfill
\begin{minipage}{0.47\linewidth}
  \centerline{\includegraphics[width=1\linewidth]{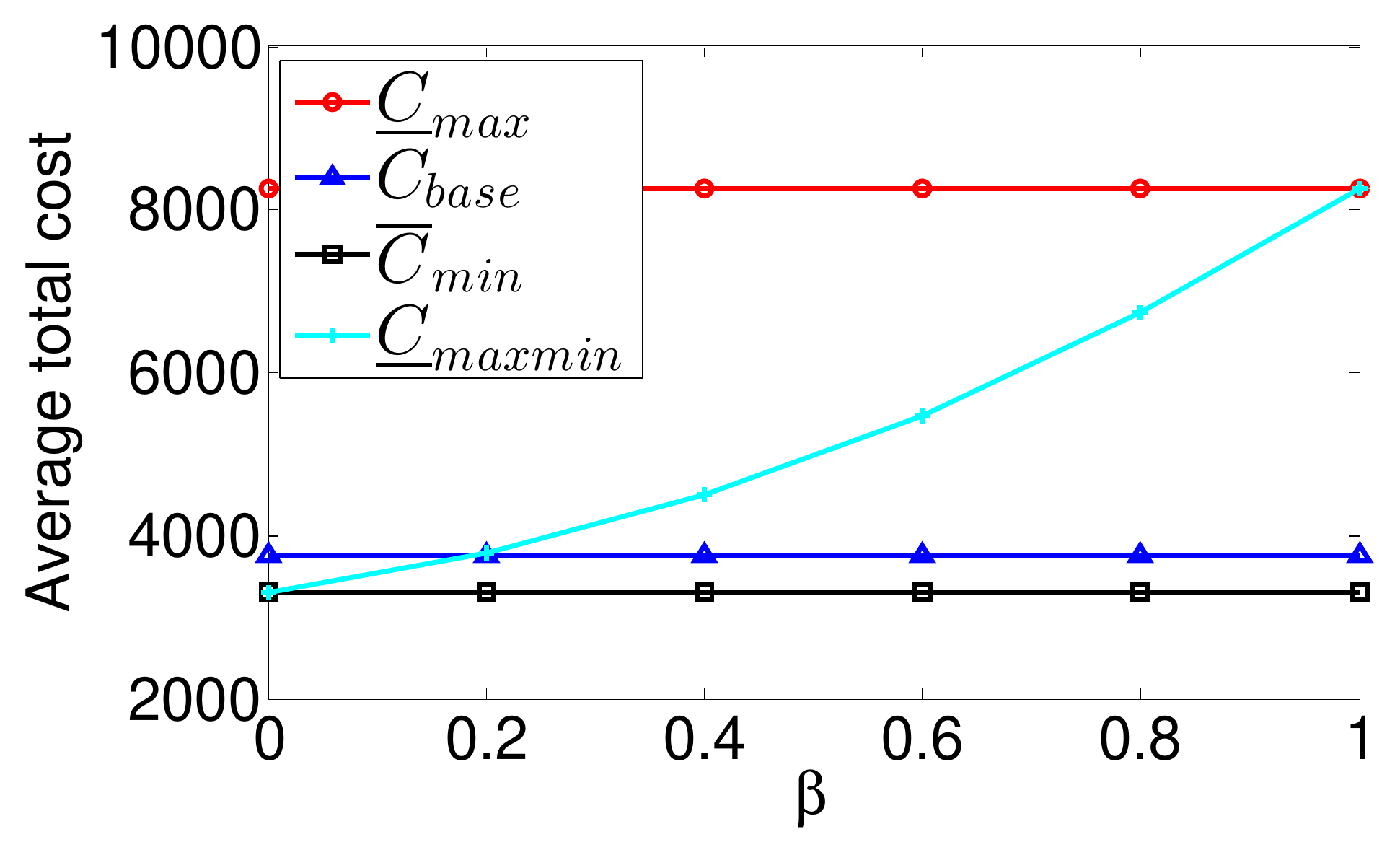}}
  \centerline{\footnotesize $x_i \in [0,40]$}
\end{minipage}\\
\begin{minipage}{0.47\linewidth}
  \centerline{\includegraphics[width=1\linewidth]{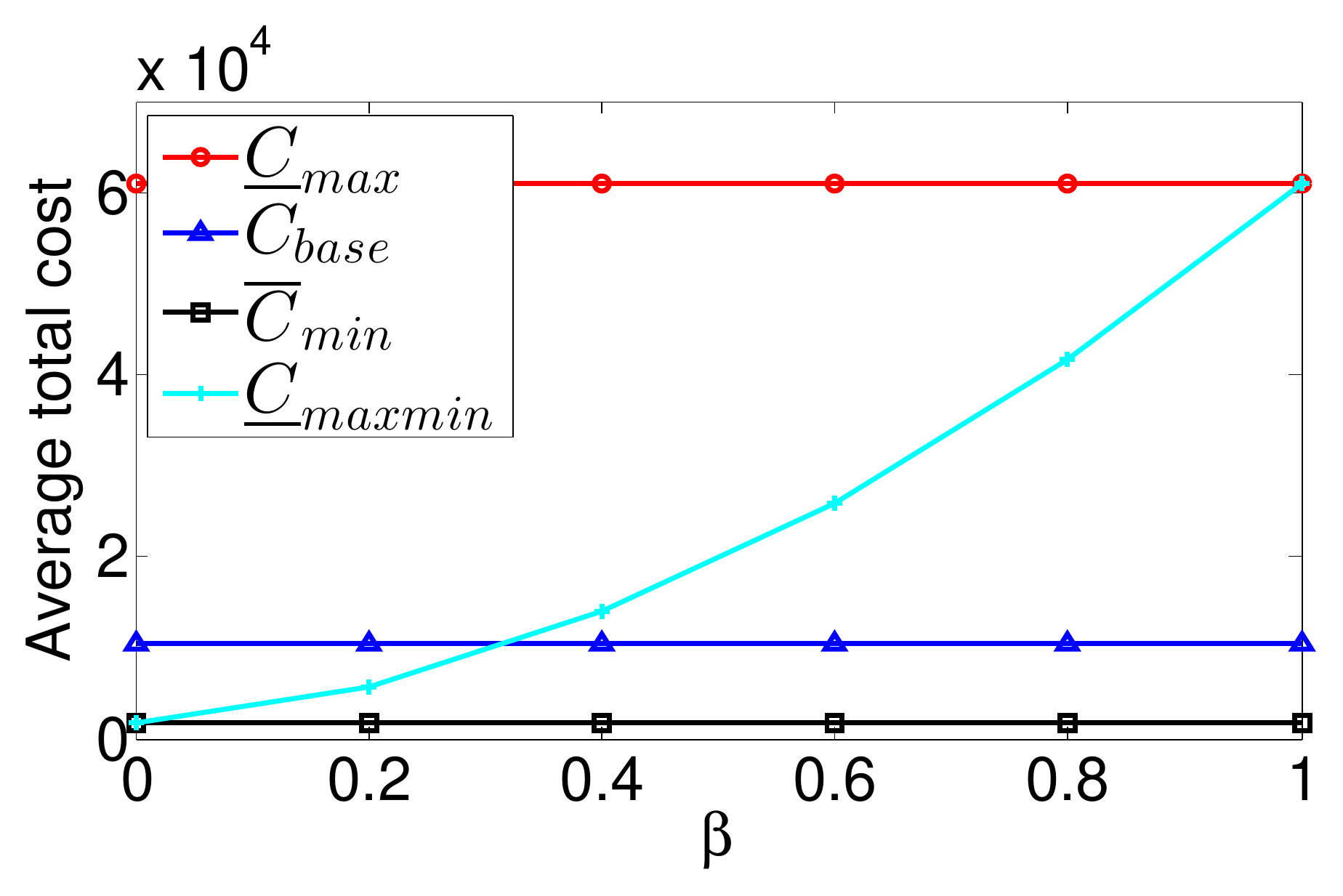}}
  \centerline{\footnotesize $x_i \in [0,10] (10\%), [40,50] (90\%)$}
  \centerline{\small (a) total-energy model}
\end{minipage}
%\hfill
\begin{minipage}{0.47\linewidth}
  \centerline{\includegraphics[width=1\linewidth]{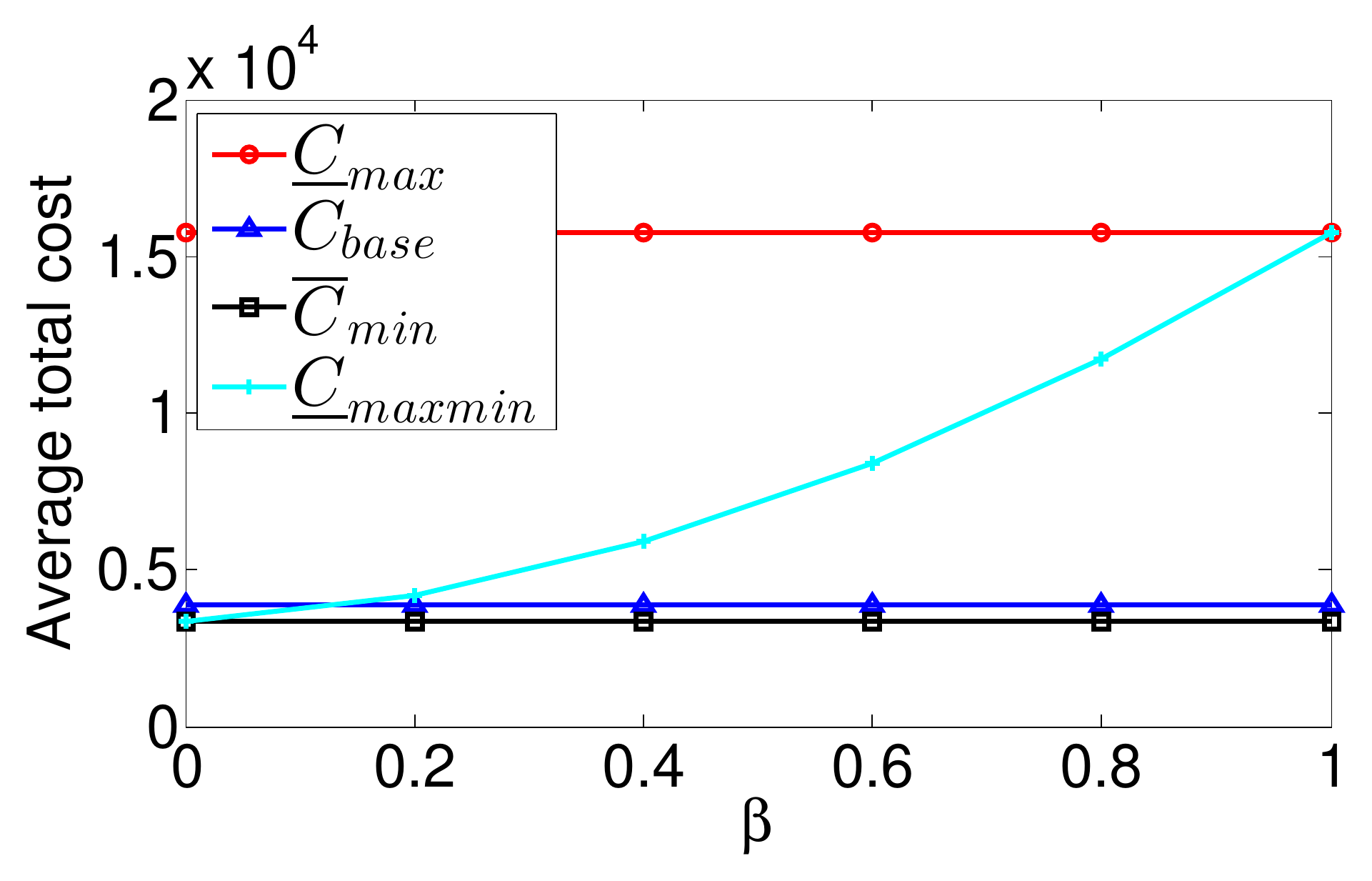}}
  \centerline{\footnotesize $x_i \in [0,10] (10\%), [40,50] (90\%)$}
  \centerline{\small (b) constant-power model}
\end{minipage}
\end{tabular}
\caption{Performance of a partially-compromised smart grid under online limited attacks with various values of $\beta$.} %The dashed lines in both figures represent the cost, $\underline{C}_{maxmin}(\beta,\delta)$, attained under online limited attacks.}
\label{fig:partial}
\end{figure}

%\noindent{\bf Offline Limited Attack:}

\vspace{1ex}
{
\noindent{\bf Online Limited Attacks:}
We now investigate the performance of online limited attacks and compare them with online full attacks. We assume that the operator schedules the set of (partially) modified demands using the AVR algorithm for the total-energy model, and the CR algorithm for the constant-power model. Since online attacks have lower complexity than their offline counterparts, we consider a larger setting with 100 jobs and each simulation is repeated 100 times. We consider the same power requirement, service time, and inter-arrival time distributions as before. Theorem~\ref{prop:max_online_tightness} and Theorem~\ref{prop:max_lower_bound} together indicate that a higher cost can be expected if most jobs have large job slackness. To confirm this, we consider two job slackness distributions, (1) a uniform distribution between [0,40], and (2) a mixture of two types of demands, where $90\%$ of demands have high elasticity with their slackness uniformly distributed in [40,50], and $10\%$ of demands are more emergent with their slackness uniformly distributed in [0,10]. The attacks were conducted with $\beta$ values ranging between 0 and 1. Figure~\ref{fig:partial} reports our results for these attacks, where the values for the corresponding online full attacks and the baselines are also plotted for reference. We observe that large job slackness can indeed enforce higher cost. For both models, even with a low fraction of jobs to be modified, the attacker still causes significant harm, compared to the un-compromised system. Moreover, the attacker becomes capable of driving the system to perform worse than its nominal point (the regular grid) with $\beta$ as low as 0.4 and 0.2 for the total-energy model and the constant-power model, respectively.
}

\iftp
\begin{figure}[!t]
\centering
\begin{tabular}{cc}
\begin{minipage}{0.47\linewidth}
  \centerline{\includegraphics[width=1\linewidth]{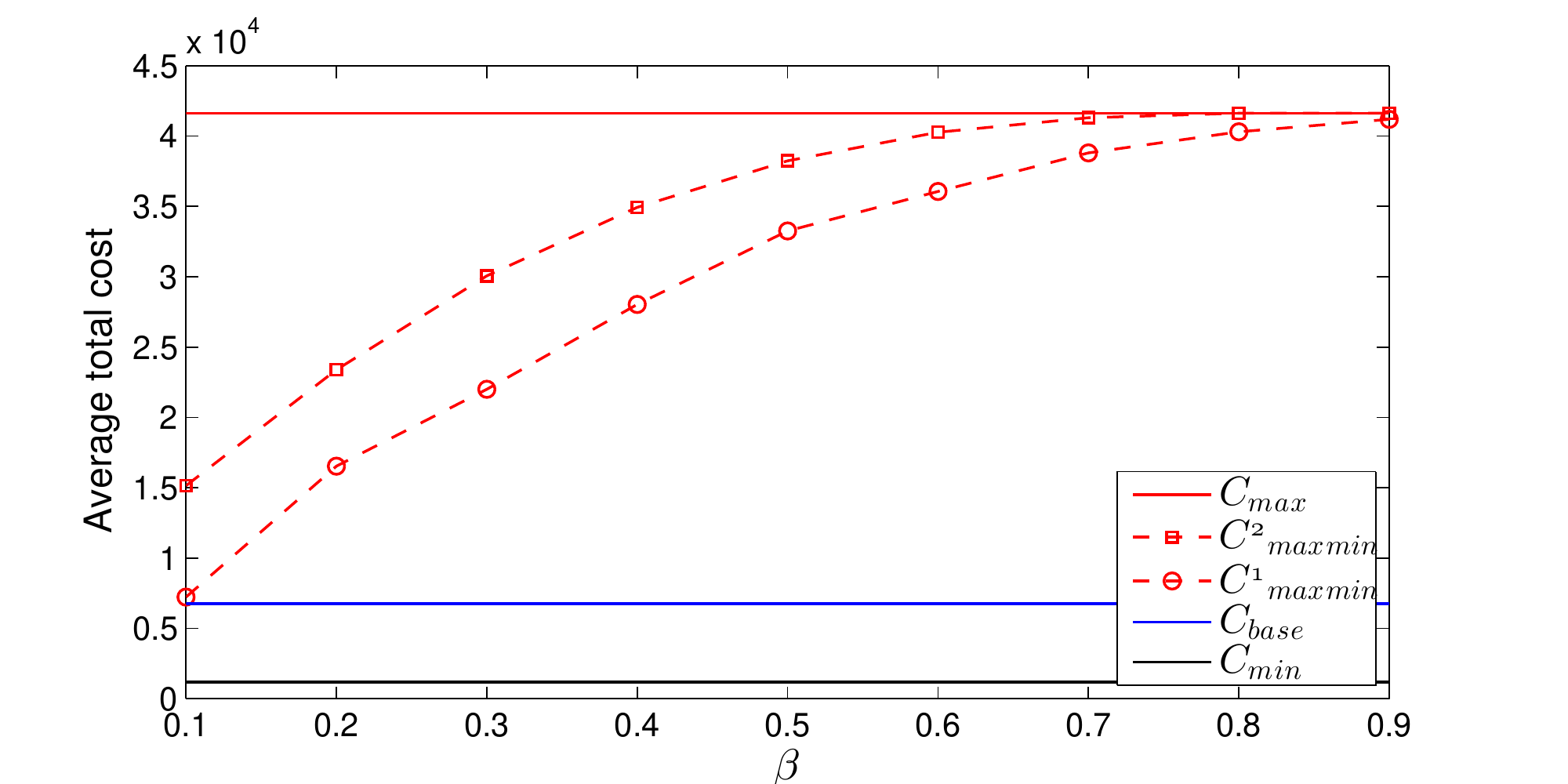}}
  \centerline{\small (a)}
\end{minipage}
%\hfill
\begin{minipage}{0.47\linewidth}
  \centerline{\includegraphics[width=1\linewidth]{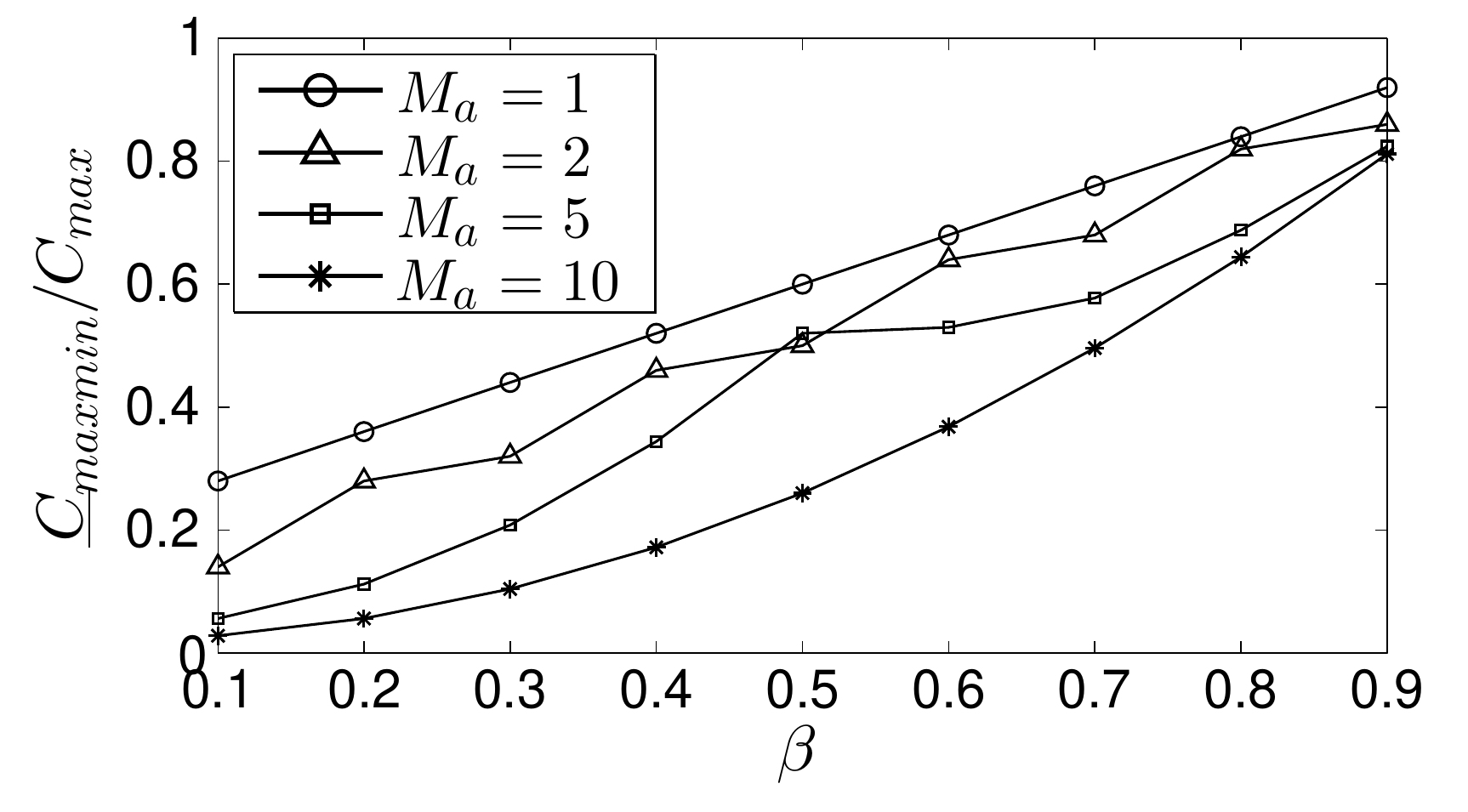}}
  \centerline{\small (b)}
\end{minipage}
\end{tabular}
\caption{Performance of offline limited attacks in the total-energy model with a varying $\beta$, for 50 jobs. In (a), energy demands are uniformly distributed on $[1,20]$, the mean interarrival time is 5, and the mean job allowance is 40. In (b), 50 identical jobs with $p_j = 5, l_j = 50$ are generated. The interarrival times are all set to $M_a$.}
\label{fig:offline-limited-TE}
\end{figure}

\vspace{0.5ex}
\noindent{\bf Offline Limited Attack:} Figure~\ref{fig:offline-limited-TE}(a) sheds more focus on the performance bounds of offline limited attacks in the total-energy model, where $C^1_{maxmin}$ and $C^2_{maxmin}$ denote the lower bound and the upper bound derived in Section~\ref{sec:te_offline_limited}, respectively. The simulation sample is composed of 50 jobs. The energy requirements were uniformly distributed on $[1,20]$ while the mean job allowance was set to 40. The results are averaged over 5 trials. As shown, with the increased allowance mean, the obtained clique partitions become denser and therefore the upper and lower bounds become tighter. Also, observe that using a simple greedy algorithm, the attacker is immediately capable of achieving a cost arbitrarily close to $C_{base}$ for our sample, with a chance of altering only 5 jobs out of 50.

Finally, we study offline limited attacks in the total-energy model in a more controlled experiment. We generate 50 identical demands, with each requiring a 5 energy units and an allowance of 50. The job interarrival times are all set to one value, denoted by $M_a$, which varies between 1 and 10. Figure~\ref{fig:offline-limited-TE}(b) shows $\underline{C}_{maxmin}/C_{max}$ under varying values of $\beta$. This enables us to gain more insights on the growth of $C_{maxmin}$ with respect to $\beta$, and how this growth is affected by the clique densities. As shown in the figure, when $M_a=1$, with our chosen parameters, a single clique of jobs could be formed to achieve the maximum cost, and hence, in accordance with our theoretical results, the attacker could achieve approximately $\beta^2$ of the maximum achievable cost. As $M_a$ increases, the growth of $\underline{C}_{maxmin}/C_{max}$ with $\beta$ approaches a linear trend. The reason is that as $M_a$ increases, the size of the optimal clique partition of jobs increases, having approximately equally sized cliques. Hence the maximum cost decreases so does the contribution of each clique to the maximum cost.
\fi

{
\section{Observations and Suggestions}~\label{sec:rec}
From our analytical studies and simulation results, we make several observations and suggestions to the operator for thwarting the new type of attacks that we consider in the paper.

\vspace{0.5ex}
\noindent{\bf Information Hiding:} We observe that the attacker's capability is significantly constrained by the amount of information it has regarding the operator and the demand patterns. In particular, to derive the best $\beta$, the attacker needs to know the intrusion detection algorithm and the key parameters such as the significance level used by the operator. Moreover, the attacker requires some prior information about the demands to make best use of its budget, such as the number of demands and the ranges of their values. Therefore, one efficient approach to reduce the damage is to properly hide these information from the attacker, e.g., by introducing noise into the data and algorithms.

\vspace{0.5ex}
\noindent{\bf Intrusion Detection:} We suggest to develop robust intrusion detection schemes that can strike a balance between the potential loss from attacks and the cost of detection. In particular, we suggest to develop a better statistical modeling of time-elastic demands, and study advanced stream data mining algorithms that can deal with the high dimension of the demand data set. Moreover, as we discussed above, it is useful to develop intrusion detection algorithms that can make it hard for the attacker to derive efficient parameters to use. %in addition to  This is in contrast to traditional intrusion detection schemes that focus on high detection accuracy.

\vspace{0.5ex}
\noindent{\bf Load Management:} We note that the scheduling algorithm used by the operator has a big impact on the total energy cost, especially when the attacker can only compromise a small number of demands. We have provided efficient solutions for the operator in the total-energy model, but better solutions are needed for the constant-power model and more general demand models. For instance, Figure~\ref{fig:partial} indicates that online limited attacks are more efficient in the constant-power model. We believe that this is due in part to the poor performance of the CR algorithm in our setting.  Moreover, it is important to develop {\it robust} algorithms that can provide a guaranteed performance even when part of demands have been modified by adversaries.

\vspace{0.5ex}
\noindent{\bf Robust and Adaptive Defense:} We suggest to develop robust defense algorithms to identify the set of most critical channels (or smart meters) to protect. From our analysis, it is clear that those demands (or a set of overlapping demands) with highest power requirement and maximum time elasticity are most beneficial to the attacker, due to the large gap between $C_{max}$ and $C_{min}$ if we consider these demands only. When these demands are mostly generated by a given subset of customers, the corresponding links can be protected to efficiently reduce damage. In the face of more advanced attackers, however, a fixed defense strategy is insufficient, as the attacker can always identify the weakest link in the system. Therefore, it is important to study adaptive defense strategies in the face of strategic attackers.
} 
\section{Conclusion}\label{sec:conclusion}
In this paper, we have studied the performance of the smart grid, in terms of energy efficiency, in the presence of an active attacks on the system. In the presence of a limited intrusion detection mechanism at the grid operator, we have proposed optimal scheduling and undetectable attack strategies. We have derived lower and upper bounds on the maximum achievable cost by an attacker with low complexity, online algorithms. %In addition, we gave bounds on the impact of attacks that are limited by intrusion detection at the operator. In these limited attacks, we have shown that a significant increase in cost could still be achieved by a simple greedy algorithm.
Overall, our theoretical analysis and numerical results show that the time-elasticity of electric load, when exploited by malicious attacks, %an inelastic utilization of the communication channels in the smart grid
could result in costs significantly higher than those expected for both the smart grid and the current electric grid, motivating the need for stronger intrusion detection and defense strategies for grid operators.

\bibliographystyle{IEEEtran}
\bibliography{refs}

\iftp\appendix
\subsection{Discussion of the Model}\label{modeling}
%As a first attempt towards understanding the impact of stealthy attacks on smart-grid demand-response, we have made several simplifications in this work. In the following, we discuss the rationale behind our model and outline several extensions.

%\vspace{1ex}
%\blue{
%\noindent{\bf Data integrity attack and intrusion detection}: In this paper, we consider a new type of data integrity attack against AMI. %On the other hand, the impact of data attacks towards AMI networks have not received much attention yet. Existing works mainly focus on cryptography based prevention techniques instead of intrusion detection.
%In contrast, we consider the problem from a decision-theoretic perspective. Although statistical testing and decision theory approach for intrusion detection in traditional networks can also be applied to AMI networks, the main challenge is that, statistics of timing can be  We have considered a simple intrusion detection scheme in this work. Moreover, all these works focus on the static setting for both the attacker and defender instead of the more realistic online setting as we consider in the paper.
%}

\vspace{1ex}
\noindent{\bf Demand-response scheme}: Our model is built upon the optimization framework proposed in~\cite{Koutsopoulos2012}. Similar models where customers submit their total energy demands together with their time elasticity have also been adopted in some recent works on electric vehicle charging~\cite{EV-LangTong-2012, EV-JAIR-2013}. %Our model departs from pricing based demand-response schemes in the literature, where customers can trade electricity usage with price.
We choose this model for the following reasons. First, various studies indicate that customers often prefer simpler pricing schemes, e.g., flat-rate pricing. Requiring every customer to submit a bidding curve as in more advanced pricing schemes may be difficult to apply in practice. Second, current pricing based demand-response schemes cannot model the time elasticity of electric load explicitly, which, however, can be utilized to reduce electricity cost and eventually benefit both the operator and the customers even under flat-rate pricing. %Therefore, we envision that our model and the insights obtained can be useful in practice.
It is an interesting problem to extend our studies to more sophisticated demand-response schemes where customers are more actively involved.

\vspace{1ex}
\noindent{\bf Forecast at the operator}: The demand/load forecast capability of the system operator could further limit stealthy attacks, which is not considered in the current model. In the extreme case when the operator knows everything about the future load, %which, in our model, includes all the electricity requirements and their time elasticity,
an attacker cannot modify any demand without of being detected. In practice, however, the system operator only has a rough estimate about future load distribution, which leaves room to stealthy attacks. %Our current model has considered some simple constraints on the attacker, while load distribution can be viewed a more advanced constraint.
It is an interesting problem to properly model the forecast capability of the operator for time-elastic electric load, and extend our framework to design stealthy attacks that can maximize energy cost while ensuring the forged demands to be still consistent with the load forecast.

{
\vspace{1ex}
\noindent{\bf Capacity constraint}: In our current model, we put no limit on the total energy served in each time slot to study the worst-case damage that a stealthy attacker can possible cause. This is also practical when there is always sufficient energy supply and the available capacities of distribution lines or transformers exceed the peak load. When the system is under congestion, however, both the operator and the attacker face more challenging optimization problems, especially in the online setting. In fact, when there is zero information on future arrivals, the only solution, if there is one, that can ensure all the demands are served by their deadlines is the Earliest Deadline First (EDF) policy, where jobs with earliest deadlines are served as fast as possible subject to the capacity constraint. To obtain a more useful problem formulation in this new setting, one approach is to relax the deadline constraints of jobs, and introduce a utility function for customers, as we further elaborate below.

\vspace{1ex}
\noindent{\bf Beyond energy cost}: We have considered two demand models with different levels of flexibility in this work. It is possible to consider more general demand models as in~\cite{lijunchen-smartgrid}, where for each customer, there is an upper and a lower bound on the energy served in each time slot, together with a utility function defined over the resulting service vector. Alternatively, we can also relax the deadline constraints by introducing a penalty for unsatisfied demands when the system is congested. A reasonable objective for the system operator is then to maximize the welfare, in terms of the total customer utility minus the total energy cost. Such flexibility provides further opportunity for the operator to improve the energy efficiency, which, however, may also be exploited by malicious attackers to harm both the system and the customers. It is interesting to extend our stealthy attack algorithms to study the fundamental tradeoffs involved in this more general setting.
}

\iftp
\subsection{Time-Dependent Cost Functions}\label{sec:time-varing-cost}
A time-invariant energy cost curve has been assumed in Section~\ref{sec:te} and Section~\ref{sec:st}. Due to the dynamics on both demand and supply, especially the uncertainty introduced by the penetration of renewable energy, energy cost can exhibit significant time variations. It is therefore important to study the impact of time-dependent cost functions on both the operator and the attacker. %Let $C_t(\cdot)$ denote the cost function at time $t$, which is assumed to be strictly convex and monotone. Let $C'_t(E)$ denote the derivative (i.e., the marginal cost) at energy load $E$.
In this section, we show that most of our previous results can be readily extended to strictly convex and monotone cost function $C_t(\cdot)$ that can vary over time. %Below we focus on the total energy model for simplicity. %As an example, we will consider $C_t(E) = c_tE^b, b \in \mathbb{R}, b \geq 1$ and derive performance bounds for it.

\subsubsection{Scheduling at the Operator} We first note that the offline YDS algorithm can be extended to time-dependent cost functions by replacing the notion of energy intensity introduced in Section~\ref{sec:te_min} by {\it energy derivative} defined below. For simplicity, we further assume that $C_t(\cdot)$ has continuous derivative, $C_t(0) = 0$, and $C'_t(0) = 0$, for any $t$. For the received (forged) demands $J'$ and an time interval $[k,l]$, let $\mathcal{S}$ denote a {\it locally} optimal schedule with minimum energy cost for jobs entirely contained in $[k,l]$, which is unique by our assumptions on $C_t(\dot)$. We define the energy derivative of the interval to be
\begin{equation}\label{eq:intensity}
    \gamma(\mathcal{I}_{J'}(k,l)) = \min_{t \in [k,l]} C'(E_\mathcal{S}(t)).
\end{equation}

That is, the energy derivative is defined as the {\it minimum} marginal cost of any time slot in $[k,l]$ in the locally optimal schedule. A critical interval is then defined as an interval with the maximum energy derivative. We observe that from our assumptions about $C_t(\cdot)$, %there is a unique optimal local schedule for a critical interval,
each time slot in a critical interval must have the same marginal cost. Moreover, by a similar argument as in~\cite{Yao1995}, it can be shown that there is an optimal schedule for all the jobs, where jobs in a critical interval is scheduled exactly as its locally optimal schedule. It follows that Algorithm~\ref{alg:min_offline} can be extended to get an optimal offline schedule for time-dependent cost by replacing energy intensity by energy derivative.

We further note that when $C_t(E) = c_tE^b, b \in \mathbb{R}, b \geq 1$, there is an online scheduling algorithm for the total energy model that achieves a competitive ratio of $O(b^b)$~\cite{Gupta-waoa2012}, assuming that upon the arrival of any job $j$, the cost functions up to $d_j$ are known to the operator. The algorithm extends AVR and looks for a minimum cost allocation for each new request on its arrival, given the previous scheduled requests while ignoring the future arrivals. %For the constant-power model, the continuous relaxation of the operator's problem can still be solved in polynomial time for time varying cost. On the other hand, whether the CR policy in~\cite{Koutsopoulos2012} can be extended to this case requires further investigation.

\subsubsection{Scheduling at the Attacker} For time-dependent cost, we show that Algorithm~\ref{alg:max_offline} can be readily extended to obtain an optimal offline attack under the total energy model. First, the optimal attack still corresponds to a clique partition of the set of jobs since Lemma~\ref{lem:compress} is proved for the general case and Lemma~\ref{lem:clique_partition} only depends on Lemma~\ref{lem:clique_partition}, although in this new setting, all the jobs in a clique should be compressed to a time slot that achieves the maximum cost among all the time slots where the job intervals intersect. For any clique $K$, let $t_K$ denote such a time slot. Define the marginal cost of $K$ as the derivative of the cost function at $t_K$ after serving all the jobs in $K$. Lemma~\ref{lem:locally_maximal} can then be proved for time-dependent cost by considering the clique $K_i$ with the maximum marginal cost in a clique partition. If $K_i$ is not locally maximal at $t_{K_i}$, then a job $j$ that intersects $t_{K_i}$ can be moved from another clique to $K_i$ without decreasing the total cost. It also follows that Theorem~\ref{thm:recursion} still holds. Hence, Algorithm~\ref{alg:max_offline} can be extended to time-dependent cost functions.

Since Algorithm~\ref{alg:max_online} is derived from Algorithm~\ref{alg:max_offline}, it can also be extended to derive online full attacks under time-dependent cost and the total energy model. Moreover, the performance bound in Theorem~\ref{prop:max_online_tightness} can be generalized as follows. Assume $C_t(E) = c_tE^b, b \in \mathbb{R}, b \geq 1$. Define $c_{max} = \max_t c_t$ and $c_{min} = \min_t c_t$. Then the online algorithm achieves at least a fraction $\frac{c_{min}}{c_{max}}\frac{1}{r^{b-1}}$ of the offline optimal cost. %Similar bounds corresponding to Theorems~\ref{prop:max_lower_bound} and~\ref{prop:max_upper_bound} can also be derived.

For limited attacks, we remark that Algorithm~\ref{alg:maxmin_greedy} can be readily extended to time-dependent cost while achieving the same lower bound as~\eqref{eq:te_maxmin_bound_1}. On the other hand, the upper bound does not apply anymore as Theorem~\ref{prop:maxmin_DP} does not hold in this new setting. For online limited attacks, Algorithm~\ref{alg:maxmin_online} can be extended to time-dependent cost. %with a performance bound similar to~\eqref{eq:te_maxmin_bound} with an extra factor of $\frac{c_{min}}{c_{max}}$, similar to the full attack case discussed above.
Finally, similar results can be derived for the constant-power model as well.
\fi

\iftp
\subsection{Algorithm 2 (Offline Full Attacks)}\label{alg_max_offline}
%\begin{alg}[\cite{Gijswijt2007}]
For all $k\in [1,T]$, set the initial condition
  \begin{equation}
     \overline{C}(k,k) = C \bigg( \sum_{j\in \mathcal{I}_J(k,k)} e_j \bigg).
  \end{equation}

With increasing interval width, iterate over all intervals $[k,l], k \leq l, k,l \in [0,T]$, and apply the following dynamic program:
\begin{enumerate}
  \item Compute
  \vspace{-2ex}
  \begin{equation*}
    \overline{C}(k,l) = \max_{z\in[k,l]} \Bigg[ C \bigg( \sum_{j\in K^z_{k,l}} e_j \bigg) + \overline{C}(k,z-1) + \overline{C}(z+1,l) \Bigg]
  \end{equation*}
   \noindent with $z^*$ achieving the optimality.

  \item Update the clique partition
  \vspace{-1ex}
  \begin{equation*}%\label{eq:rec02}
    \mathcal{Q}(k,l) = \begin{cases}
    \emptyset, \qquad \mbox{if } \mathcal{I}_J(k,l) = \emptyset,\\
    \mathcal{Q}(k,z^*-1) \cup K^{z^*}_{k,l} \cup \mathcal{Q}(z^*+1,l), \hspace{1ex} \mbox{otherwise.}
    \end{cases}
  \end{equation*}
\end{enumerate}
%\end{alg}
\fi

\iftp
\else
\begin{figure}
 \centering
 \includegraphics[width=0.4\textwidth]{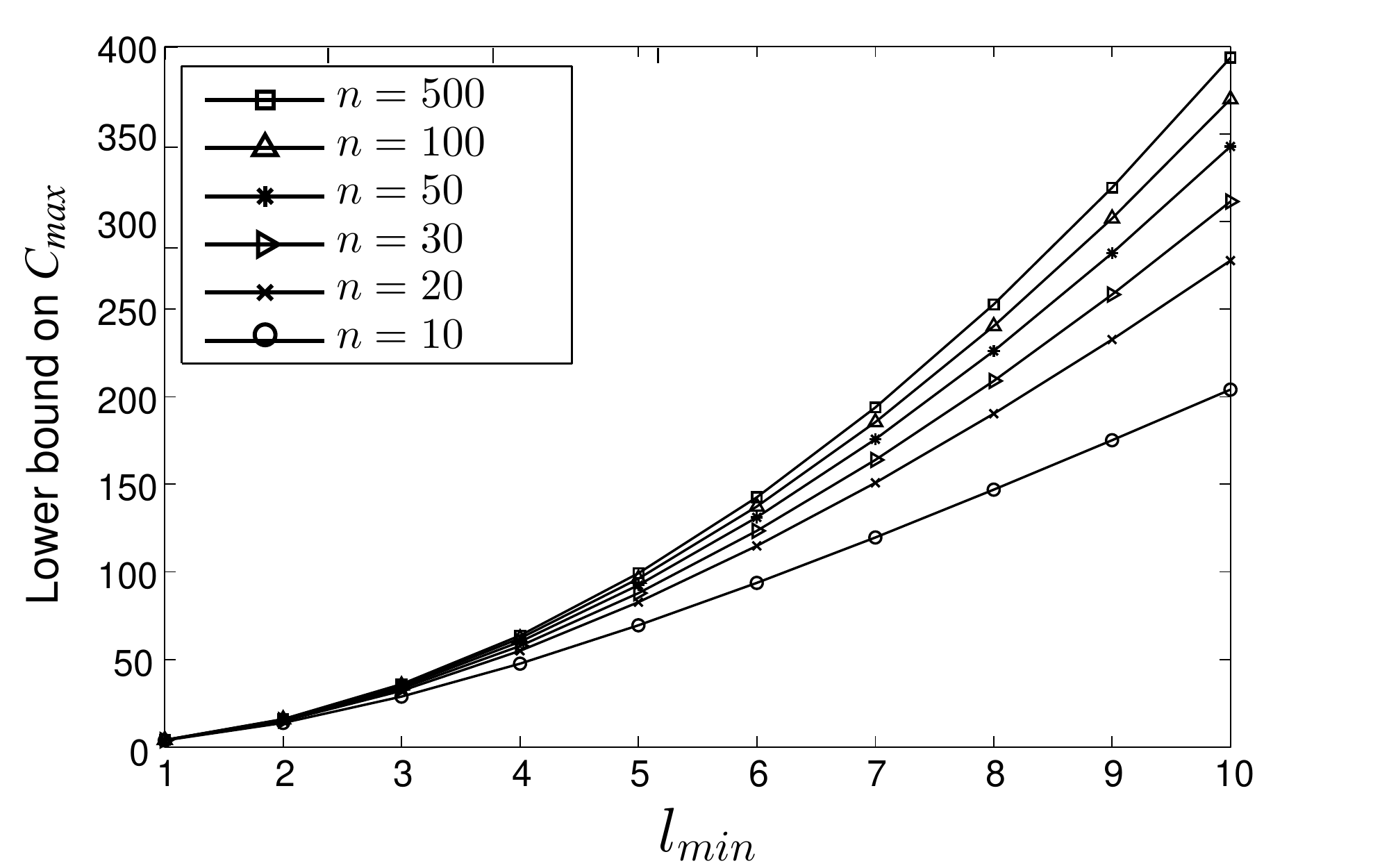}
 \caption{A lower bound on $C_{max}$ plotted for various values of $n$ and $l_{min}$ under a quadratic cost function (i.e., $b = 2$). The average energy demand is 10 while the average inter-arrival time is 5.}
 \label{fig:max_lower_bound}
\end{figure}
\fi

\fi
\end{document}